\documentclass[12pt]{article}
\def\today{November 29, 2003}
\usepackage{amsmath}
\usepackage{amssymb}
\usepackage{amsthm}
\usepackage{psfig}
\newtheorem{thm}{Theorem}[section]
\newtheorem{co}[thm]{Corollary}
\newtheorem{lem}[thm]{Lemma}

\newtheorem{assumption}[thm]{Assumption}

\newtheorem{definition}[thm]{Definition}
\newenvironment{de}{\begin{definition}\rm}{\end{definition}}
\newtheorem{example}[thm]{Example}
\newenvironment{exmp}{\begin{example}\rm}{\end{example}}
\newtheorem{remark}[thm]{Remark}
\newenvironment{rem}{\begin{remark}\rm}{\end{remark}}
\newtheorem{tab}{Table}
\newenvironment{ta}{\begin{tab}\rm}{\end{tab}}
\newcommand{\length}{{\rm length}\,}

\newcommand{\vier}[4]{\left( \begin{array}{ccc}
                   #1 &\;& #2 \\ #3 &\;& #4 \end{array} \right)}

\newcommand{\Section}[1]{\section{#1}\setcounter{equation}{0}}

\newcommand{\eqr}[1]{~\mbox{$(${\rm \ref{#1}}$)$}}
\newcommand{\tr}{{\rm tr}\,}
\newcommand{\diag}{{\rm diag}\,}

\newcommand{\V}{\mathcal{V}}
\renewcommand{\S}{\mathcal{S}}
\newcommand{\R}{\mathbb{R}}
\newcommand{\C}{\mathbb{C}}

\title{Geometrical and Numerical Design of Structured \\
Unitary Space Time Constellations \footnote{Both authors were
supported in part by NSF grants DMS-00-72383 and CCR-02-05310. The
first author was also supported by a fellowship from the Center of
Applied Mathematics at the University of Notre Dame. A preliminary
version of this paper was presented at 40-th Allerton Conference
on Communication, Control, and Computing, Monticello, Illinois,
October 2002.}}

\date{{\normalsize \today}}

\author{Guangyue Han,\ \ Joachim Rosenthal\\
  {\normalsize Department of Mathematics}\vspace{-1mm} \\
  {\normalsize University of Notre Dame}\vspace{-1mm} \\
  {\normalsize Notre Dame, IN 46556.}\\
  {\normalsize {\em e-mail:\/} Han.13@nd.edu,\  Rosenthal.1@nd.edu}\vspace{-1mm}\\
  {\normalsize {\em URL:} http://www.nd.edu/\~{}eecoding/} }

\begin{document}\maketitle\thispagestyle{empty}

%%%%%%%%%%%%%%%%%%%%%%%%%%%%%%%%%%%%%%%%%%%%%%%%%%%%%%%%%%%%
\begin{abstract}
  Unitary space-time modulation using multiple antennas promises
  reliable communication at high transmission rates. The basic
  principles are well understood and certain criteria for
  designing good unitary constellations have been presented.

  There exist two important design criteria for unitary space
  time codes. In the situation where the signal to noise ratio is
  large it is well known that the {\em diversity product} (DP) of
  a constellation should be as large as possible. It is less
  known that the {\em diversity sum} (DS) is a very important
  design criterion for codes working in a low SNR environment.
  For some special situations, it will be more practical and
  reasonable to consider a constellation optimized at a certain SNR
  interval. For this reason we introduce the {\em diversity
    function} as a general design criterion.  So far, no general
  method to design good-performing constellations with large
  diversity for any number of transmit antennas and any
  transmission rate exists.

  In this paper we propose constellations with suitable structure which
  allow one to construct codes with excellent diversity using
  geometrical symmetry and numerical methods. We also demonstrate
  how these structured constellations out-perform currently existing
  constellations and explain why the proposed constellation structure admit
  simple decoding algorithm: sphere decoding. The presented design
  methods work for any dimensional constellation and for any
  transmission rate. Moreover codes based on the proposed structure
  are very flexible and can be optimized for any signal to noise
  ratio.
\end{abstract}

\newpage
%%%%%%%%%%%%%%%%%%%%%%%%%%%%%%%%%%%%%%%%%%%%%%%%%%%%%%%%%%%%
\Section{Introduction and Model}

One way to acquire reliable transmission with high transmission
rate on a wireless channel is to use multiple transmit or receive
antennas. Either because of rapid changes in the channel
parameters or because of limited system resources, it is
reasonable to assume that both the transmitter and the receiver
don't know about the channel state information (CSI), i.e. the
channel is non-coherent.

In~\cite{ho00a}, Hochwald and Marzetta study unitary space-time
modulation. Consider a wireless communication system with $M$
transmit antennas and $N$ receive antennas operating in a
Rayleigh flat-fading channel. We assume time is discrete and at
each time slot, signals are transmitted simultaneously from the
$M$ transmitter antennas. We can further assume that the wireless
channel is quasi-static over a time block of length $T$.

A signal constellation $\V:=\{ \Phi_1,\ldots, \Phi_L\}$ consists
of $L$ matrices having size $T \times M$ and satisfying $T \ge M$
and $\Phi_k^* \Phi_k = I_M$. The last equation simply states that
the columns of $\Phi_k$ form a ``unitary frame'', i.e. the column
vectors all have unit length in the complex vector space
$\mathbb{C}^T$ and the vectors are pairwise orthogonal.  The
scaled matrices $\sqrt{T} \Phi_k$, $k=1,2,\cdots,L$, represent
the code words used during the transmission. It is known that the
transmission rate is determined by $L$ and $T$:

$$
\mathtt{R}=\frac{\log_2(L)}{T}.
$$

Let $\rho$ represent the expected signal-to-noise ratio (SNR) at
each receive antenna. The basic equation between the received
signal $R$ and the transmitted signal $\sqrt{T} \Phi$ is given
through:
$$
R=\sqrt{\frac{\rho T}{M}}\Phi H+W,
$$
where the $M \times N$ matrix $H$ accounts for the
multiplicative complex Gaussian fading coefficients and the $T
\times N$ matrix $W$ accounts for the additive white Gaussian
noise. The entries $h_{m,n}$ of the matrix $H$ as well as the
entries $w_{t,n}$ of the matrix $W$ are assumed to have a
statistically independent normal distribution
$\mathcal{CN}(0,1)$. In particular it is assumed that the
receiver does not know the exact values of either the entries of
$H$ or $W$ (other than their statistical distribution).

The decoding task asks for the computation of the most likely
sent code word $\Phi$ given the received signal $R$. Denote by
$||\ \ ||_F$ the Frobenius norm of a matrix. If $A=(a_{i,j})$
then the Frobenius norm is defined through $|| A
||_F=\sqrt{\sum_{i,j} |a_{i,j}|^2}.$ Under the assumption of the
above model the maximum likelihood (ML) decoder will have to
compute:
$$
\Phi_{ML}=\displaystyle \arg \max_{\Phi_l \in
  \{\Phi_1,\Phi_2,\cdots,\Phi_L\}} {\|R^*\Phi_l\|}_F
$$
for each received signal $R$. (See~\cite{ho00a}).

Let $\delta_m(\Phi_l^* \Phi_{l'})$ be the $m$-th singular value
of $\Phi_l^* \Phi_{l'}$. It has been shown in~\cite{ho00a} that
the pairwise probability of mistaking $\Phi_l$ for $\Phi_{l'}$
using maximum likelihood decoding satisfies:

\begin{eqnarray}
P_{\Phi_l,\Phi_{l'}} &=& \mbox{Prob}\left(\mbox{ choose }\Phi_{l'}\mid
\Phi_{l}\mbox{ transmitted } \right)(\rho)\nonumber\\
 &=& \mbox{Prob}\left(\mbox{ choose }\Phi_{l}\mid
\Phi_{l'}\mbox{ transmitted } \right)(\rho) \nonumber\\
&=& \frac{1}{4\pi} \int_{-\infty}^{\infty}\frac{4}{4w^2+1}\prod_{m=1}^M
\left[1+
\frac{(\rho T/M)^2 (1-\delta_m^2 (\Phi_l^*
\Phi_{l'})) }{4(1+\rho T/M)} (4w^2+1)\right]^{-N}\!\! dw \label{exactf}\\
 &\le& \frac{1}{2} \prod_{m=1}^M
\left[1+
\frac{(\rho T/M)^2(1-\delta_m^2 (\Phi_l^*
\Phi_{l'}))}{4(1+\rho T/M)} \right]^{-N}. \label{mainf}
\end{eqnarray}
It is a basic design objective to construct constellations $\V=\{
\Phi_1,\ldots, \Phi_L\}$ such that the pairwise probabilities
$P_{\Phi_l,\Phi_{l'}}$ are as small as possible. Mathematically
we are dealing with an optimization problem with unitary
constraints:

\vspace{0.3cm}

Minimize $\displaystyle \max_{l \neq l'} P_{\Phi_l,\Phi_{l'}}$
with the constraints $\Phi_i^*\Phi_i=I$ where $i=1,2,\cdots,L$.

\vspace{0.3cm}

Formula\eqr{mainf} is sometimes referred to as ``Chernoff's
bound''.  This formula is easy to work with, the exact
formula\eqr{exactf} is in general not easy to work with, although
it could be useful in the numerical search of good constellations
as well. Researchers have been searching for constructions where
the maximal pairwise probability of $P_{\Phi_l,\Phi_{l'}}$ is as
small as possible. Of course the pairwise probabilities depend on
the chosen signal to noise ratio $\rho$ and the construction of
constellations has therefore to be optimized for particular
values of the SNR.

The design objective is slightly simplified if one assumes that
transmission operates at high SNR situations. In~\cite{ho00}, a
design criterion for high SNR is presented and the problem has
been converted to the design of a finite set of unitary matrices
whose diversity product is as large as possible. In this special
situation several researchers~\cite{al98,ta00a,sh01,sh02} came up
with algebraic constructions and we will say more about this in
the next section.

The main purpose of this paper is to present structured
constellation and to develop geometrical and numerical procedures
which allow one to construct unitary constellations with excellent
diversity for any set of parameters $M,N,T,L$ and for any signal
to noise ratio $\rho$. The paper is structured as follows. In
Section~\ref{Sec-diversity} we introduce the diversity function
of a constellation. This function depends on the signal to noise
ratio and it gives for each value $\rho$ an indication how well
the constellation $\V$ will perform. For large values of $\rho$
the diversity function is governed by the diversity product, for
small values of $\rho$ it is governed by the diversity sum. These
concepts are introduced in Section~\ref{Sec-diversity} as well.
The introduced concepts are illustrated on some well known
constellations previously studied in the literature.

In Section~\ref{Sect-Alg} we first show that randomly constructed
codes are fully diverse with probability one. Then we start the
main task of this paper, namely to parameterize constellations
which will be efficient for numerical search algorithms. For this
purpose we introduce the concept of a {\em weak group structure}
and we classify all weak group structures whose elements are
normal and positive.

In Section~\ref{Sec-geometrical} we investigate an algebraic
structure which led to some of the best constellations which we
were able to derive. We also show that in the good-performing
codes the distance spectrum profile for both the diversity sum
and the diversity product are important.

Section~\ref{Sec-numerical} is one of the main sections of this
paper.  We first explain a general method on how one can
efficiently design excellent constellations for any set of
parameters $M,N,T,L$ and $\rho$. For this we review the
properties of the complex Stiefel manifold and the Cayley
transform. We conclude this section with an
extensive table where we publish a large set of codes having some
of the best diversity sums and diversity products in their
parameter range. More extensive lists of codes with large
diversity can be found on the website~\cite{ha03u2}.

Finally in Section~\ref{sphere-decoding} we explain how the
algebraic structure which underlies most of the derived codes can
be used to have a fast decoding algorithm.  Our simulations
indicate that in the design of codes more attention should be
given to the diversity sum (more generally diversity function)
which previously has not been fully studied.

%%%%%%%%%%%%%%%%%%%%%%%%%%%%%%%%%%%%%%%%%%%%%%%%%%%%%%%%%%%%%%%%%%%%%
\Section{The Diversity Function, the Diversity Product (DP) and
the Diversity Sum (DS)} \label{Sec-diversity}

In this paper we will be concerned with the construction of
constellations where the right hand sides in\eqr{exactf}
and\eqr{mainf}, maximized over all pairs $l,l'$ is as small as
possible for fixed numbers of $T,M,N,L$. As already mentioned this
tasks depends on the signal to noise ratio the system is
operating.  For this purpose we define the {\em exact
  diversity function} dependent on the constellation $\V=\{
\Phi_1,\ldots, \Phi_L\}$ and a particular SNR $\rho$ through:
\begin{equation}                  \label{exact-div}
\mathcal{D}_e(\V,\rho):=
\max_{l \ne l'} \mbox{Prob}\left(\mbox{ choose }\Phi_{l'}\mid
\Phi_{l}\mbox{ transmitted } \right)(\rho)
\end{equation}
For a particular constellation with a large number $L$ of
elements, with many transmit and receive antennas the function
$\mathcal{D}_e(\V,\rho)$ is very difficult to compute. Indeed for
each pair $\Phi_{l'},\Phi_{l}$ it is required to compute the
singular values of the $M\times M$ matrix $\Phi_l^* \Phi_{l'}$ and
then one has to evaluate up to $L(L-1)/2$ integrals of the
form\eqr{exactf} and this has to be done for each value of $\rho$.
Although this task is formidable it can be done in cases where
$T,M,L$ are all in the single digits using e.g. Maple.

\begin{figure}[ht]
\centerline{\psfig{figure=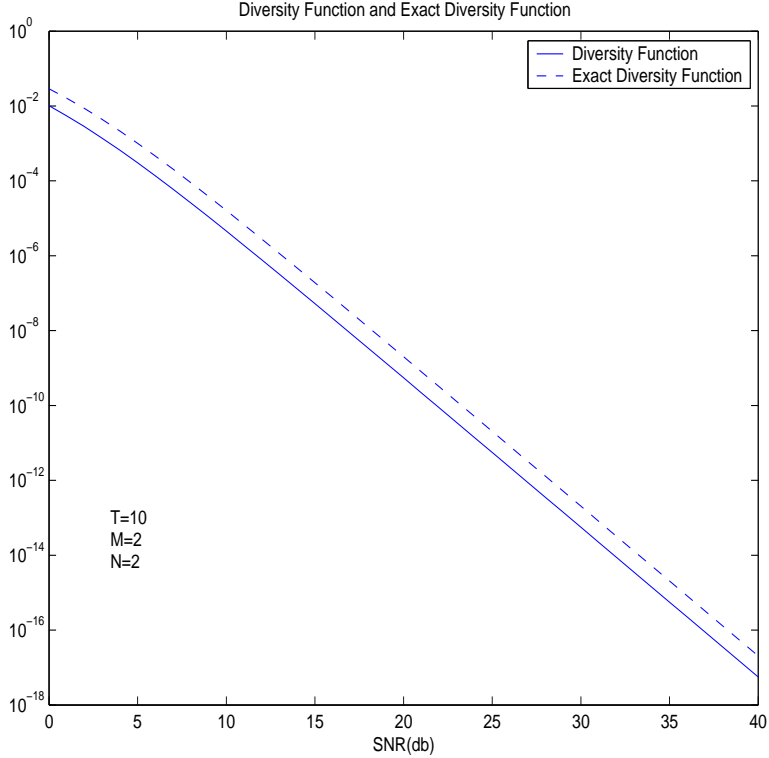,width=4in,height=4in}}
\caption{Diversity function $\mathcal{D}(\V,\rho)$ and exact
  diversity function $\mathcal{D}_e(\V,\rho)$ of a fully
  orthonormal constellation.} \label{fig-2}
\end{figure}

Using Chernoff's bound\eqr{mainf} we define a simplified function
called the {\em diversity function} through:
\begin{equation}                  \label{div}
\mathcal{D}(\V,\rho):=
\max_{l \ne l'}
\frac{1}{2} \prod_{m=1}^M
\left[1+
\frac{(\rho T/M)^2}{4(1+\rho T/M)} (1-\delta_m^2 (\Phi_l^*
\Phi_{l'}))\right]^{-N}.
\end{equation}
The computation of $\mathcal{D}(\V,\rho)$ does not require the
evaluation of an integral and the computation requires
essentially the computation of $ML(L-1)/2$ singular values. The
singular values $\delta_m (\Phi_l^*\Phi_{l'})$ are by definition
all real numbers in the interval $[0,1]$ as we assume that the
columns of $\Phi_l,\Phi_{l'}$ form both orthonormal frames.  The
functions $\mathcal{D}_e(\V,\rho)$ and $\mathcal{D}(\V,\rho)$ are
the smallest if the singular values $\delta_m (\Phi_l^*\Phi_{l'})$
are as small as possible. These numbers are all equal to zero if
and only if the column spaces of $\Phi_l,\Phi_{l'}$ are pairwise
perpendicular. We call such a constellation {\em fully
orthonormal}. Since the columns of $\Phi_l$ generate an
$M$-dimensional subspace this can only happen if $L\leq T/M$. On
the other hand if $L\leq T/M$ it is easy to construct a
constellation where the singular values of $(\Phi_l^*\Phi_{l'})$
are all zero. Just pick $LM$ different columns from a $T\times T$
unitary matrix. Figure~\ref{fig-2} depicts the functions
$\mathcal{D}_e(\V,\rho)$ and $\mathcal{D}(\V,\rho)$ for a fully
orthonormal constellation with $T=10$ and $M=N=2$.

In order to study the function $\mathcal{D}(\V,\rho)$ more
carefully let
\begin{equation}                              \label{tilderho}
\tilde{\rho}:=\frac{(\rho T/M)^2}{4(1+\rho T/M)}.
\end{equation}
In some small interval $[\rho_1,\rho_2]$ the maximum
in\eqr{div} is achieved for some fixed indices $l,l'$ and in
terms of $\tilde{\rho}$ the function $\mathcal{D}(\V,\rho)$ is of
the form:
$$
\mathcal{D}(\V,\tilde{\rho})=\frac{1}{2\left(
    1+c_1\tilde{\rho}+\cdots +c_M\tilde{\rho}^M\right)^N},
$$
where the coefficients $c_1,\ldots,c_M$ depend on the
particular constellation and on the chosen interval
$[\rho_1,\rho_2]$. For an interval close to zero the dominating
term will be the coefficient~$c_1$. Up to some factor this term
will define the {\em diversity sum} of the constellation. When
$\tilde{\rho}>>0$ then the dominating term will be the
coefficient $c_M$ and up to some scaling this term will define
the {\em diversity product} of the constellation. A constellation
will have a small diversity function for small values of $\rho$
(and presumably performs well in this range) when the
constellation is chosen having a large diversity sum.  A
constellation will have a small diversity function for large
values of $\rho$ (and presumably performs well in this range)
when the constellation is chosen having a large diversity
product.  In the next two subsections we will study the limiting
behavior of $\mathcal{D}(\V,\rho)$ as $\rho$ goes to zero and to
infinity.

%%%%%%%%%%%%%%%%%%%%%%%%%%%%%%%%%%%%%%%%%%%%%%%%%%%%%%%%%%%%%%%%%%%%%
\subsection{Design criterion for high SNR}

When the SNR $\rho$ is very large then $\mathcal{D}(\V,\rho)$ can
be approximated via:
\begin{equation}
\mathcal{D}(\V,\rho)\simeq
\max_{l \ne l'}
\frac{1}{2}
\left(
\frac{(\rho T/M)^2}{4(1+\rho T/M)}
\right)^{-NM}
\prod_{m=1}^M
\frac{1}{(1-\delta_m^2 (\Phi_l^*
\Phi_{l'}))^{N}}.
\end{equation}
It is the design objective to construct a constellation $\Phi_1,
\Phi_2,\cdots, \Phi_n$ such that
$$
\min_{l \ne l'} \prod_{m=1}^M(1-\delta_m^2(\Phi_l^*
\Phi_{l'}))
$$
is as large as possible.  This last expression defines in
essence the diversity product.  In order to compare different
dimensional constellations it is customary to use the definition:
\begin{de} (See~\cite{ho00})                             \label{div-prod}
  The {\em diversity product} of a unitary constellation $\V$ is
  defined as
  $$
  \prod \V = \min_{l \ne l'} \left(\prod_{m=1}^M (1-
    \delta_m(\Phi_l^* \Phi_{l'}) ^2)\right)^{\frac{1}{2M}}.
  $$
\end{de}

An important special case occurs when $T=2M$. In this situation
it is customary to represent all unitary matrices $\Phi_k$ in the
form:
\begin{equation}                         \label{specialform}
\Phi_k=\frac{\sqrt{2}}{2} \left(\begin{array}{c}
    I\\
    \Psi_k
  \end{array}\right).
\end{equation}
Note that by definition of $\Phi_k$ the matrix $\Psi_k$ is a $M
\times M$ unitary matrix. The diversity product as defined in
Definition~\ref{div-prod} has then a nice form in terms of the
unitary matrices. For this let $\lambda_m$ be the $m$th
eigenvalue of a matrix, then
$$
1-\delta_m^2(\Phi_{l'}^* \Phi_l)=\frac{1}{4}
\lambda_m(2I_M-\Phi_l^* \Phi_{l'}-\Phi_{l'}^* \Phi_l)
=\frac{1}{4}\delta_m^2(I_M-\Psi_{l'}^*
\Psi_l)=\frac{1}{4}\delta_m^2(\Psi_{l'}-\Psi_l).
$$
So we have
$$
\prod_{m=1}^M (1-\delta_m^2(\Phi_{l'}^*
\Phi_l))^{\frac{1}{2M}}=\frac{1}{2}\prod_{m=1}^M
\delta_m(\Psi_{l'}-\Psi_l)^{\frac{1}{M}}=\frac{1}{2}|
\det(\Psi_{l'}-\Psi_l)|^{\frac{1}{M}}.
$$
When $T=2M$ and the constellation $\V$ is defined as above,
then the formula of the diversity product assumes the simple
form:
\begin{equation}  \label{diversity}
\prod \V=\frac{1}{2} \min_{0 \leq l < l' \leq L}
|\det(\Psi_l-\Psi_{l'})|^{\frac{1}{M}}.
\end{equation}

We call a constellation $\V$ a fully diverse constellation if
$\prod \V > 0$. A lot of efforts have been taken to construct
constellations with large diversity product. (See
e.g.~\cite{ho00,li02,ha02p2,ha02p,sh01,sh02,ta00a}). For the
particular situation $T=2M$ with special form\eqr{specialform}
the design asks for the construction of a discrete subset
$\V=\{\Psi_1,\ldots,\Psi_L\}$ of the set of $M\times M$ unitary
matrices $U(M)$. When this discrete subset has the structure of a
discrete subgroup of $U(M)$ then the condition that $\V$ is fully
diverse is equivalent to the condition that the identity matrix is
the only element of $\V$ having an eigenvalue of 1. In other
words the constellation $\V$ is required to operate fixed point
free on the vector space $\mathbb{C}^M$. Using a classical
classification result of fixed point free unitary representations
by Zassenhaus~\cite{za36}, Shokrollahi et al.~\cite{sh01,sh02}
were able to study the complete list of fully diverse finite group
constellations inside the unitary group $U(M)$. Some of these
constellations have the best known diversity product for given
fixed parameters $M,N,L$. Unfortunately the possible
configurations derived in this way is somehow limited. The
constellations are also optimized for the diversity product and
as we demonstrate in this paper for unitary space time modulation
maybe attention should be given to the diversity sum.

In most of the literature mentioned above researchers focus their
attention to constellations having the special
form\eqr{specialform}. Unitary differential modulation~\cite{ho00}
is used to avoid sending the identity (upper part of every element
in the constellation) redundantly. This increases the
transmission rate by a factor of 2 to:

$$
\mathtt{R}=\frac{\log_2(L)}{M}=2\frac{\log_2(L)}{T}.
$$

\noindent Because of this reason we will also focus ourselves in
the later part of the paper to the special form\eqr{specialform}
as well. Nonetheless it will become obvious that the numerical
techniques also work in the general situation.

%%%%%%%%%%%%%%%%%%%%%%%%%%%%%%%%%%%%%%%%%%%%%%%%%%%%%%%%%%%%
\subsection{Design criterion for low SNR channel}

As we mentioned before a constellation with a large diversity sum
will have a small diversity function at small values of the
signal to noise ratio. This is particularly suitable when the
system operates in a very noisy environment. When $\rho$ is
small, using Formula\eqr{tilderho}, one has the following
expansion:
\begin{multline*}
\prod_{m=1}^M [1+ \frac{(\rho T/M)^2}{4(1+\rho T/M)} (1-\delta_m^2
(\Phi_l^* \Phi_{l'}))]= \prod_{m=1}^M [1+ \tilde{\rho}(1-
\delta_m^2(\Phi_l^* \Phi_{l'}))]\\
=1+\tilde{\rho} \sum_{m=1}^M (1-\delta_m^2(\Phi_l^*
\Phi_{l'}))+O(\tilde{\rho}^2).
\end{multline*}

When $\rho \rightarrow 0$, i.e. $\tilde{\rho} \rightarrow 0$, we
can omit the higher order terms $O(\tilde{\rho}^2)$ and the upper
bound of $P_{\Phi_l,\Phi_{l'}}$ requires that
$$
\sum_m (1-\delta_m^2(\Phi_l^* \Phi_{l'}))=(M-{\|\Phi_l^* \Phi_{l'}
\|}_F^2)
$$
is large.  In order to lower the pairwise error probability, it is
the objective to make ${\|\Phi_l^* \Phi_{l'} \|}_F^2$ as small as
possible for every pair of $l, l'$. It follows that at high SNR,
the probability primarily depends on $\prod_{m=1}^M
(1-\delta_m^2(\Phi_l^* \Phi_{l'}))$, but at low SNR, the
probability primarily depends on $\sum_{m=1}^M
(1-\delta_m^2(\Phi_l^* \Phi_{l'}))$. In order to be able to
compare the constellation of different dimensions, we define:
\begin{de}                             \label{div-sum}
The {\em diversity sum} of a unitary constellation $\V$ is defined
as
$$
\sum \V = \min_{l \ne l'} \sqrt{1-\frac{
{\|\Phi_l^*\Phi_{l'}\|}_F^2}{M}}.
$$
\end{de}

Again one has the important special case where $T=2M$ and the
matrices $\Phi_k$ take the special form\eqr{specialform}. In this
case one verifies that
\begin{multline*}
{\|\Phi_{l}^*\Phi_{l'}\|}_F^2=\frac{1}{4}{\|I+\Psi_l^*\Psi_{l'}\|}_F^2
=\frac{1}{4}\tr((I+\Psi_{l'}^*\Psi_l)(I+\Psi_l^*\Psi_{l'}))\\
=\frac{1}{4}\tr(2I+\Psi_{l'}^*\Psi_l+\Psi_l^*\Psi_{l'})=
\frac{1}{4}(4M-(2M-\tr(\Psi_{l'}^*\Psi_l+\Psi_l^*\Psi_{l'})))\\
=\frac{1}{4}(4M-\tr((\Psi_l-\Psi_{l'})^*(\Psi_l-\Psi_{l'})))=
\frac{1}{4}(4M-{\|\Psi_l-\Psi_{l'}\|}_F^2)
\end{multline*}

For the form\eqr{specialform} the diversity sum assumes the
following simple form:
\begin{equation}                                         \label{T2M-div-sum}
\sum \V = \min_{l,l'}\frac{1}{2\sqrt{M}}{\|\Psi_l- \Psi_{l'} \|_F}.
\end{equation}
Without mentioning the term the concept of diversity sum was used
in~\cite{ho00a1}.  Liang and Xia~\cite[p. 2295]{li02} explicitly
defined the diversity sum in the situation when $T=2M$ using
equation\eqr{T2M-div-sum}. Definition~\ref{div-sum} naturally
generalizes the definition to arbitrary constellations.

We want to point out that the diversity sum is the design criterion
only for unitary constellation. Hochwald and
Marzetta~\cite{ho00a} calculate the non-coherent space time
channel capacity and indicate that unitary signal constellation
are capacity achieving signal sets only for high SNR scenarios.
For low SNR case, the transmitting power should be allocated
unsymmetrically, i.e., unitary constellations are not capacity
achieving in the first place. However unitary signal sets are
easily manageable and one can take advantage of differential
modulation technique~\cite{ho00} to speed up the transmission.
Moreover our simulation results indicate that codes with near
optimal diversity sum tend to perform significantly better
compared to the currently existing ones optimized for the
diversity product for low and even moderate SNR scenarios.
So it is quite reasonable and more toward the practical use to construct
unitary constellations with good diversity sum.

As  the formulas make it clear the diversity sum and the
diversity product are in general very different. There is however
an exception. When $T=4$, $M=2$ and the constellation $\V$ is in
the special\eqr{specialform}. If in addition all the $2\times 2$
matrices $\{ \Psi_1,\ldots, \Psi_L\}$
are a subset of the special unitary group
$$
SU(2)=\{ A\in \C^{2\times 2}\mid A^*A=I\mbox{ and }\det A=1\}
$$
then it turns out that the diversity product $\prod \V$ and
the diversity sum $\sum \V$ of such a constellation are the same.
For this note that elements $\Psi_l,\Psi_{l'}$ of $SU(2)$ have
the special form:
$$
\Psi_l=\vier{a}{b}{-\bar{b}}{\bar{a}},\
\Psi_{l'}=\vier{c}{d}{-\bar{d}}{\bar{c}}.
$$
Through a direct calculation one verifies that
$\det(\Psi_l-\Psi_{l'})=|a-c|^2+|b-d|^2$ and
${\|\Psi_l-\Psi_{l'}\|}_F^2=2(|a-c|^2+|b-d|^2)$. But this means
that $\prod \V= \sum \V$ for constellations inside $SU(2)$.

%%%%%%%%%%%%%%%%%%%%%%%%%%%%%%%%%%%%%%%%%%%%%%%%%%%%%%%%%%%%%%%%%%
\subsection{Four illustrative examples}    \label{SubSec-I}

The diversity sum and the diversity product govern the diversity
function at low SNR respectively at high SNR. Codes optimized at
these extreme values of the SNR-axis do not necessarily perform
well on the ``other side of the spectrum''. In this subsection we
illustrate the introduced concepts on four examples. All examples
have about equal parameters, namely $T=4$, $M=2$ and the size $L$
is 121 respectively 120. The first two examples are well studied
examples from the literature. We derived the third and the fourth
examples by geometrical design and numerical methods respectively.

\paragraph{Orthogonal Design:} This constellation has been
considered by several authors~\cite{al98,sh01}. For our purpose
we simply define this code as a subset of $SU(2)$:
$$
\left\{\frac{\sqrt{2}}{2}\left(\begin{array}{cc}
            e^{\frac{2m\pi i}{11}}&e^{\frac{2n\pi i}{11}}\\
            -e^{-\frac{2n\pi i}{11}}&e^{-\frac{2m\pi i}{11}}
            \end{array} \right)|m,n=0,1,\cdots,10\right\}.
$$
The constellation has 121 elements and the diversity sum and the
diversity product are both equal to $0.1992$.

\begin{figure}[ht]
\centerline{\psfig{figure=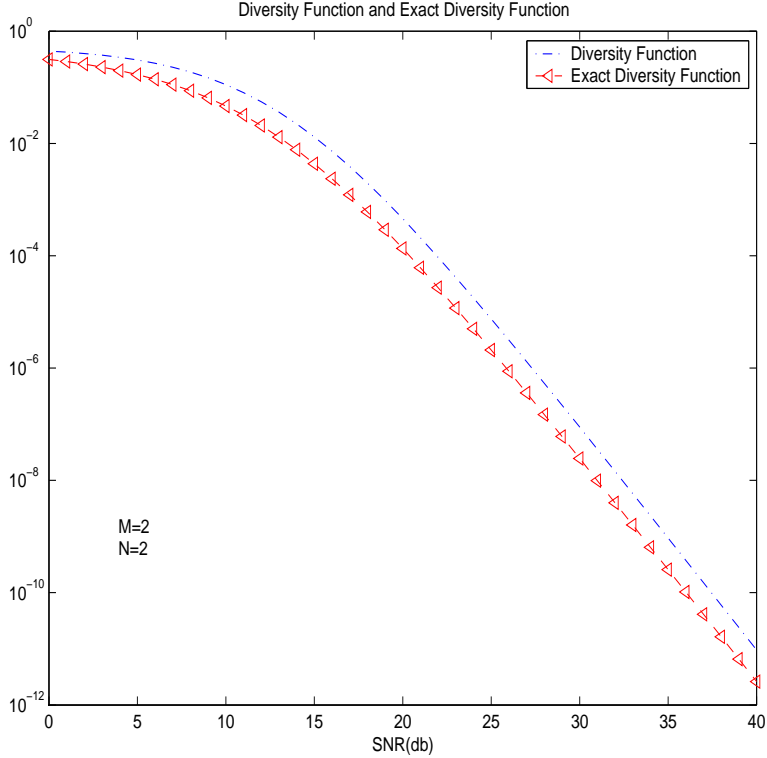,width=4in,height=4in}}
\caption{Diversity function $\mathcal{D}(\V,\rho)$ and exact
  diversity function for the group constellation
$SL_2(\mathbb{F}_5)$.}\label{fig-1}
\end{figure}

\paragraph{Unitary Representation of $SL_2(\mathbb{F}_5)$:}
Shokrollahi et al.~\cite{sh01} derived a constellation using the
theory of fixed point free representations whose diversity
product is near optimal. This constellation appears as a unitary
representation of the finite group $SL_2(\mathbb{F}_5)$ and we
will refer to this constellation as the
$SL_2(\mathbb{F}_5)$-constellation. The finite group
$SL_2(\mathbb{F}_5)$ has 120 elements and this is also the size
of the constellation. In order to describe the constellation let
$\eta = e^{\frac{2\pi i}{5}}$ and define
$$
P=\frac{1}{\sqrt{5}} \left(\begin{array}{cc}
    \eta^2-\eta^3&\eta^1-\eta^4\\
    \eta^1-\eta^4&\eta^3-\eta^2\\
 \end{array} \right),\ \
Q=\frac{1}{\sqrt{5}} \left(\begin{array}{cc}
    \eta^1-\eta^2&\eta^2-\eta^1\\
    \eta^1-\eta^3&\eta^4-\eta^3\\
\end{array} \right).
$$
Then the constellation is given by the set of matrices
$(PQ)^jX$,where $j = 0,1,\cdots,9$, and $X$ runs over the set
\begin{multline*}
\{ I_2, P, Q, QP, QPQ, QPQP, QPQ^2, QPQPQ, QPQPQ^2, \\
QPQPQ^2P, QPQPQ^2PQ, QPQPQ^2PQP \}.
\end{multline*}
The constellation has rate $R = 3.45$ and $\prod
{SL_2(\mathbb{F}_5)}=\sum {SL_2(\mathbb{F}_5)} =
\frac{1}{2}\sqrt{\frac{(3-\sqrt{5})}{2}} \sim 0.3090$. The
diversity product of this constellation is truly outstanding. For
illustrative purposes we plotted in Figure~\ref{fig-1}  the
exact diversity functions and the diversity function of this constellation.

\paragraph{Numerically Derived Constellation:}
Using simulated annealing algorithm we found after short
computation a constellation with very good diversity sum. The
constellation is given through a set of 121 matrices
\begin{multline*}
\left\{\Psi_{k,l}:=A^kB^l|A=\left(\begin{array}{cc}
                 -0.9049 + 0.3265*i&   0.1635 + 0.2188*i\\
                 0.0364 + 0.2707*i&  -0.8748 + 0.4002*i
                 \end{array}\right),\right. \\
\left.
B=\left(\begin{array}{cc}
        -0.1596 + 0.9767*i&  -0.1038 + 0.0994*i\\
        0.0833 - 0.1171*i&  -0.9432 + 0.2995*i
        \end{array}\right), k,l=0,1,\cdots,10 \right\}.
\end{multline*}
As we explain in Section~\ref{sphere-decoding}, the maximum
likelihood decoding of this constellation admits a simple decoding
algorithm: sphere decoding.
\bigskip

\paragraph{Geometrically Designed Constellation:}
Based on the algebraic structure we are going to propose in this
paper, we further implement the geometrical symmetry into this
structure. A geometrically designed constellation can be
described as follows:
\begin{multline*}
\left\{\Psi_{k}:=A^kB^k|A=\left(\begin{array}{cc}
                 e^{17 \pi/60}&   0\\
                 0&  e^{13 \pi/60}
                 \end{array}\right),\right. \\
\left. B=\left(\begin{array}{cc}
        \cos(22\pi/60)&  \sin(22\pi/60)\\
        -\sin(22\pi/60)&  \cos(22\pi/60)
        \end{array}\right), k=0,1,\cdots,119 \right\}.
\end{multline*}
This constellation has superb diversity sum and reasonably good
diversity product. One can also use sphere decoding to implement maximum likelihood
decoding of this constellation.
\bigskip

\begin{ta}
The following table summarizes the parameters of the four
constellations:
\begin{center}
\begin{tabular}{|c|c|c|c|c|}
  \hline
   & \begin{tabular}{c}
Orthogonal\\design\end{tabular}
 & $SL_2(\mathbb{F}_5)$
& \begin{tabular}{c} Numerically\\ derived\end{tabular} &
\begin{tabular}{c}
Geometrically\\ designed\end{tabular} \\
  \hline
  Number of elements& 121 & 120 & 121 & 120\\
  \hline
  diversity sum&  0.1992 &0.309 & 0.3886 & 0.4156  \\
  \hline
  diversity product & 0.1992&0.309 &  0.0278 & 0.1464 \\ \hline
\end{tabular}
\end{center}
\end{ta}

\begin{figure}[ht]
\centerline{\psfig{figure=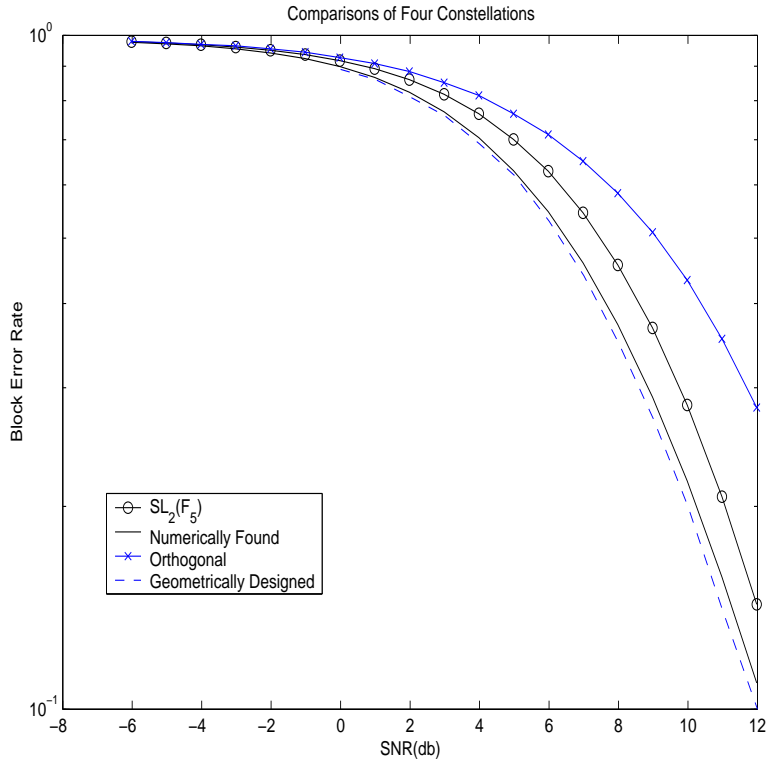,width=4in,height=4in}}
\caption{Simulations of four constellations having sizes
  $T=4$, $M=2$ and $L=120$ respectively $L=121$.}
\label{fig-3}
\end{figure}

Of course we were curious about the performances of these four
different codes. Figure~\ref{fig-3} provides simulation results
for each of the four constellations. Note that the numerically
designed code who has a very bad diversity product is performing
very well nevertheless due to the exceptional diversity sum. One
can see that up to $12$db numerically derived codes outperform
the group code by about $1$ db. In fact, our simulation
results show that until $35$db the numerical one is still
performing much better than the orthogonal one. However at around
$18$db, the group constellation surpasses the numerical one due to
exceptional diversity product. The geometrically designed
constellation has better diversity sum and diversity product than
the numerical one, therefore its performance is better than the
numerical one (our results show that their performance curves are
quite close, although the geometrical one is slightly better).
These simulation results give an indication that the diversity
sum is a very important parameter for a unitary constellation at
low SNR regime.

%%%%%%%%%%%%%%%%%%%%%%%%%%%%%%%%%%%%%%%%%%%%%%%%%%%%%%%%%%%%%%%%%%%%%%
\Section{Constellations With Algebraic Structure}
\label{Sect-Alg}

Before we venture into the realm of structured constellation, we
would like explore random unitary space time constellations first.
We introduce the Haar distributed random matrix, which in some sense
can be viewed as a high dimensional generalization of a complex
random variable with circular symmetric distribution
$\mathcal{CN}(0,1)$.

\begin{de}
  The  Haar measure on $U(M)$ is defined to be a probability measure
  $\mathcal{H}$ on $U(M)$ which is translate invariant: for any
  mensurable set $S$ in $U(M)$ and any fixed element $U_0$ in
  $U(M)$
  $$
  \mathcal{H}(S)=\mathcal{H}(U_0S).
  $$
  A unitary random matrix $\textbf{U}$ is Haar distributed
  (h.d.)  if for any measurable set $S$ we have
  $$
  Pr(\textbf{U} \in S)=\mathcal{H}(S).
  $$
\end{de}

\begin{rem}
  Note that h.d. matrix is also called isotropically distributed
  matrix in~\cite{ma99a}. We want to point out that Haar measure
  can be defined more generally. In fact every compact Lie group
  admit a unique (up to scalar) translate invariant measure: Haar
  measure~\cite{bo86b}.
\end{rem}

A well known yet non-trivial fact is that for any measurable set
$S \subset U(M)$, we have
$$
\mathcal{H}(S)=\mathcal{H}(S^*),
$$
where $S^*$ consists of the conjugate transpose of all the
elements in $S$. Thus for a h.d. matrix $\textbf{U}$, one can
verify
$$
Pr(\textbf{U}^* \in S)=Pr(U \in
S^*)=\mathcal{H}(S^*)=\mathcal{H}(S).
$$
Immediately we conclude $\textbf{U}^*$ is also h.d. matrix.
Also one can verify that the product of two h.d. matrices is still
h.d.  Another very interesting property about a h.d. matrix is
about its spectrum. As derived in~\cite{go98}, the joint
probability density for the eigenvalues of a h.d. random matrix
$\textbf{U} \sim \diag(e^{i \mathbf{\theta_1}}, e^{i
  \mathbf{\theta_2}}, \cdots, e^{i \mathbf{\theta_M}})$ in $U(M)$
is given by the Weyl denominator formula:
$$
f(\theta_1, \theta_2, \cdots, \theta_M)=\frac{1}{(2\pi)^M M!}
\prod_{j < k} {|e^{i\theta_j}-e^{i\theta_k}|}^2.
$$
The properties of h.d. matrices lead to the following theorem
about random unitary space time constellation:
\begin{thm}
  For a random unitary space time constellation $\V$ consisting
  of $L$ h.d. independent random matrices
  $\mathbf{U_1},\mathbf{U_2}, \cdots, \mathbf{U_L}$, we have
  $$
  Pr(\prod \V=0)=0,
  $$
  that is the probability of $\V$ being fully diverse is $1$.
\end{thm}

\begin{proof}
  First we can rewrite
  $$
  Pr(\prod \V=0)=Pr(\bigcup_{j < k}
  |\det(\mathbf{U_j}-\mathbf{U_k})|=0) \leq \sum_{j < k}
  Pr(|\det(\mathbf{U_j}-\mathbf{U_k})|=0).
  $$
  Next we are going to show that the probability of the event
  $|\det(\mathbf{U_j}-\mathbf{U_k})|=0$ happening is $0$. Now,
  $$
  Pr(|\det(\mathbf{U_j}-\mathbf{U_k})|=0)=
  Pr(|\det(I-\mathbf{U_j}^*\mathbf{U_k})|=0).
  $$
  Let $\mathbf{U}$ denote $\mathbf{U_j}^*\mathbf{U_k}$, we
  know $\mathbf{U}$ is h.d. matrix. Using the Weyl denominator
  formula, one computes
  $$
  Pr(|\det(I-\mathbf{U})|=0)=Pr(\bigcup_{l=1}^M
  \mathbf{\theta_l}=0) \leq \frac{1}{(2\pi)^M M!} \sum_{l=1}^M
  \int\!\!\!\int_{\theta_l=0} \prod_{j < k}
  {|e^{i\theta_j}-e^{i\theta_k}|}^2 d\theta_1 d\theta_2 \cdots
  d\theta_M.
  $$
  Since
  $$
  \int\!\!\!\int_{\theta_j=0} \prod_{j < k}
  {|e^{i\theta_j}-e^{i\theta_k}|}^2 d\theta_1 d\theta_2 \cdots
  d\theta_M \leq 2^{M(M-1)} \int\!\!\!\int_{\theta_j=0} d\theta_1
  d\theta_2 \cdots d\theta_M=0,
  $$
  we conclude that
  $$
  Pr(|\det(\mathbf{U_j}-\mathbf{U_k})|=0)=0.
  $$
  Consequently
  $$
  Pr(\prod \V=0)=0,
  $$
  that is the probability of $\V$ being fully diverse is $1$.
\end{proof}

Note that if an $M \times M$ matrix $G$ of independent complex
Gaussian entries is input to $QR$ algorithm, the resulting
unitary matrix $Q$ is Haar distributed~\cite{ea83}. For
simplicity we sketch the proof as follows: First one can write
$Q=GR^{-1}$, then for a fixed unitary matrix $U_0$, it can be
checked that $U_0 G$ has the same distribution as $G$.
Consequently $U_0 Q$ has the same distribution as $Q$, i.e., the
distribution of $Q$ is translate invariant. Therefore the
uniqueness of translate invariant measure on a compact Lie group
guarantees that $Q$ is Haar distributed. As a consequence of the
above theorem, an algorithm which produces a fully diverse unitary
constellation with probability $1$ can be given as follows: take
$L$ instance of complex Gaussian matrices and feed them through the $QR$
algorithm, the resulting $L$ unitary matrices constitute a fully
diverse constellation with probability $1$.

From an algebraic geometry point of view one easily shows that
the set of constellations with $\prod \V=0$ forms a lower
dimensional proper algebraic sub-variety of $U(M)^L$. In
particular the set of all the fully diverse constellations is
Zariski open~\cite{ha77} in $U(M)^L$, i.e., fully diverse
constellations are dense in $U(M)^L$. Haar distributed
random constellations won't be practical for maximum
likelihood decoding in high transmission rate scenario because no
algebraic structure is assumed for random constellation and
therefore the decoding process will be too complex. In the sequel
we are going to investigate structured constellations and explain
how one can restrict the parameter space to judiciously chosen
subsets and how one can convert maximum likelihood decoding to
lattice decoding by using structured constellations.

Consider a general constellation of square unitary matrices,
$$
\V=\{\Psi_1, \Psi_2, \cdots, \Psi_L\}.
$$
In order to calculate the diversity product, one needs to do
$\frac{L(L-1)}{2}$ calculations: $ |\det(\Psi_i-\Psi_j)|$ for
every different pair $i,j$. The same statement can be said about
the diversity sum, however for simplicity we only show the
diversity product case in the sequel unless specified otherwise.

If one deals with a group constellation then one needs only to
calculate $L-1$ such determinant calculations and this is one of
the remarkable advantages of group constellations. This is a
direct consequence of
$$
|\det(\Psi_i-\Psi_j)|=|\det(\Psi_i)\det(I-\Psi_i^*\Psi_j)|
=|\det(I-\Psi_i^*\Psi_j)|,
$$
where $\Psi_i^*\Psi_j$ is still in the group.

As we mentioned before group constellations are however very
restrictive about what the algebraic structure is concerned. In
the following we are going to present some constellations which
have some small number of generators and whose diversity can be
efficiently computed. This will ensure that the total
parameter space to be searched is limited as well. We start with
an example:
\begin{exmp}                           \label{Exmp4}
  Consider the constellation
  $$
  \V = \{A^kB^l|A, B \in U(M), k=0, \cdots, p, l=0, \cdots, q
  \}.
  $$
  The parameter space for this constellation is $U(M)\times
  U(M)$, this is a manifold of dimension~$2M^2$ and the number of
  elements in $\V$ is $(p+1)(q+1)$. If one has to compute
  $|\det(\Psi_i-\Psi_j)|$ for every distinct pair this would
  require $\left(\begin{array}{c} (p+1)(q+1) \\ 2
  \end{array}\right)$ determinant calculations. We will show in the following that the
same result can be obtained by doing $2pq+p+q$ determinant
computations.

Let $\Psi_i$ and $\Psi_j$ be two distinct elements having the
form $A^{k_1}B^{l_1}$ and $A^{k_2}B^{l_2}$ respectively. We have
now several cases. When $k_1 =k_2$, then necessarily $l_1 \neq
l_2$ and the distance is computed as
$$|\det(A^{k_1}B^{l_1}-A^{k_2}B^{l_2})|=|\det(I-B^{|l_2-l_1|})|,$$
where $|l_2-l_1|$ is an integer between $1$ and $q$. If $l_1
=l_2$, then we have $k_1 \neq k_2$ and the distance is computed
as
$$|\det(A^{k_1}B^{l_1}-A^{k_2}B^{l_2})|=|\det(I-A^{|k_2-k_1|})|,$$
where $|k_2-k_1|$ is an integer between $1$ and $p$. If $(k_1 <
k_2 \;\; \mbox{and} \;\; l_1 < l_2)$ or $(k_1 > k_2 \;\;
\mbox{and} \;\; l_1 > l_2)$, we have
$$|\det(A^{k_1}B^{l_1}-A^{k_2}B^{l_2})|=|\det(I-A^{|k_2-k_1|}B^{|l_2-l_1|}),$$
where $1\leq |k_2-k_1|\leq p$ and $1\leq |l_2-l_1|\leq q$.
Similarly if $(k_1 < k_2 \;\; \mbox{and} \;\; l_1 > l_2)$ or
$(k_1 > k_2 \;\; \mbox{and} \;\; l_1 < l_2)$ then
$$
|\det(A^{k_1}B^{l_1}-A^{k_2}B^{l_2})|=|\det(A^{|k_2-k_1|}-B^{|l_2-l_1|})|,
$$
with $1\leq |k_2-k_1|\leq p$ and $1\leq |l_2-l_1|\leq p$. The
total number of distances to be computed is in total equal to
$2pq+p+q$.
\end{exmp}

The number of distances to be computed indicates how complex the
calculation for the diversity is. In fact the smaller this number is,
intuitively the larger possibility of finding a unitary
constellation with good diversity we will have. An immediate
observation is that for two pair of unitary matrices $(A, B), \;
\; (C, D)$, if $(C, D)=(UAV, UBV)$ or $(C, D)=(UA^{-1}V,
UB^{-1}V)$, then $|\det(C-D)|=|\det(A-B)|$. We are going to
consider several constellations starting with this observation.

\begin{exmp}
  Consider the case that $G \subset U(n)$ is a subgroup with $L$
  elements, then for any two distinct elements $A, B \in G$ we
  have $|\det(A-B)|=|\det(I-A^{-1}B)|$ with $A^{-1}B \in G$.
  Therefore at most $L-1$ distance calculations are needed to
  derive the diversity product. The product of two group
  constellations has the similar property. Consider $G_i \subset
  U(M)$ with order $l_i$, where $i=1,2$. Let
  $$
  G=\{AB|A \in G_1, B \in G_2\}.
  $$
  Since
  $|\det(A_1B_1-A_2B_2)|=|\det(I-A_1^{-1}A_2B_2B_1^{-1})|$ with
  $A_1^{-1}A_2B_2B_1^{-1} \in G$, at most $L-1$ calculations are
  needed in this case, where $L=l_1l_2$.
\end{exmp}

\begin{exmp}
  Consider a constellation with the following form:
  $$
  \{A^iB^j|i=0,\cdots,l_1-1, j=0, \cdots, l_2-1 \;\;
  \mbox{and} \;\; A, B\in U(M), A^{l_1}=I, B^{l_2}=I\}.
  $$
  It can be checked that for the above constellation at
  most $L-1$ calculations are needed, where $L=l_1l_2$.
\end{exmp}

Group structures do have certain advantages for constructing
unitary constellations: it is less complex to calculate the
diversity product (or sum); the possibility of finding a large
diversity constellation intuitively may be increased. However
the constellations found by this approach~\cite{sh01} are really
few and far between. Somehow one wonders if the group
structure is too restrictive to find a good-performing
constellation.

In the sequel we are going to loosen the constraints imposed by
the group structures. As demonstrated in Example~\ref{Exmp4} it is
desirable to have a small dimensional manifold (in
Example~\ref{Exmp4} it was $U(M)\times U(M)$) which parameterizes
a set of potentially interesting constellations. Having such a
parameterization will help to avoid the problem of ``dimension
explosion''.  The set of constellations parameterized by
$U(M)\times U(M)$ in Example~\ref{Exmp4} are interesting as we
are not required to compute all pairwise distances in order to
compute the diversity product (sum).

\begin{de}
  Let $X$ be the set $\{x_1, x_2, \cdots, x_n\}$ and $F$ be the
  free group on the set $X$. A subset $G\subset U(M)$ is called
  {\em freely generated} if there are elements $\{g_1, g_2,
  \cdots, g_n\}\subset G$ such that the homomorphism $ \phi: F
  \longrightarrow G$ with $\phi(x_i)=g_i$ is an isomorphism.
\end{de}

An immediate consequence of this definition is that every element
in $G$ can be uniquely written as a product of $g_i$'s and
$g_i^{-1}$'s. The elements $g_i$ are called the generators of
$G$.  A freely generated subset $G$ is simply parameterized by
the set:
$$
\left\{ a_1^{p_1}a_2^{p_2}\cdots a_k^{p_k}\mid a_i \;\; \mbox{is
one of} \;\; g_i's, \; p_i\in\mathbb{Z}\right\}.
$$

Take an element $g \in G$ with its representation
$g=\prod_{i=1}^k a_i^{p_i}$, we say that the presentation is {\em
reduced} whenever $a_i \neq a_{i+1}$ for $i=1,\ldots,n-1$. Observe
that taking the product of distinct matrices $\prod_{i=1}^n A_i$
is numerically expensive, however taking the power of one matrix
$A^k$ is much easier (note that for $A=U\sum U^{-1}$ with $\sum$
diagonal, we have $A^k=U\sum^kU^{-1}$). Moreover by considering
the powers of one matrices, we are able to impose the lattice
structure to the constellation, which makes sphere decoding of
structured constellations possible. (see
Section~\ref{sphere-decoding}) Therefore we are interested in
``normal'' elements of $G$.

\begin{de}
We say that an element $g=\prod_{i=1}^k a_i^{p_i}$ in reduced
form is a {\em normal element} whenever $a_i\neq a_j$ for $i\neq
j$. A subset $\V$ of the freely generated set $G$ is said to be a
{\em normal constellation} if every non-identity element in $\V$
is normal.
\end{de}

Since finding an inverse of a matrix is numerically expensive, we
also limit our searches to positive constellations:

\begin{de}
  An element $g$ in $G$ with the reduced form $g=\prod_{i=1}^k
  {a_i}^{p_i}$ is said to be a {\em positive element} if $p_i >
  0$ for $i=1, 2, \cdots, k$.  A subset $\V$ of the freely
  generated set $G$ is said to be a {\em positive constellation}
  if every non-identity element in $\V$ is positive.
\end{de}

Positive normal constellations are desirable for numerical
searches as they can be efficiently parameterized and searched.
If one wants to compute the diversity product (or sum) of an
arbitrary positive constellation with $L$ elements one still has
to compare a total of $\binom{L}{2}$ pairs of matrices. In the
sequel we will impose more structure on a constellation
$\V\subset G$ which will guarantee that only $L-1$ pair of
elements have to be compared during the diversity product (sum)
computation.
\begin{de}
  Two unitary matrices $A, B \in G$ are said to be {\em
    equivalent} (denote by $A \sim B$) if there is unitary matrix $U\in G$ such that $A
  = UBU^{-1}$ or $A=UB^{-1}U^{-1}$. $[A]$ will denote all the
  matrices that are equivalent to $A$. For a constellation $\V
  \subset G$, we say $\V=\{ \Psi_1, \Psi_2, \cdots, \Psi_L \}$
  has a {\em weak group structure} if for any two distinct
  elements $\Psi_i, \Psi_j$ the product $\Psi_i^{-1} \Psi_j$ is
  equivalent to some $\Psi_k$.
\end{de}
The reader verifies that we indeed defined an equivalence
relation.  Note also that $\V$ has a group structure as soon as
$\Psi_i^{-1} \Psi_j$ is always another element of $\V$ and this explains
our wording.
\begin{lem}
  Let $\V=\{\Psi_0=I, \Psi_1, \Psi_2, \cdots, \Psi_{L-1}\}$ be a
  constellation with a weak group structure. In order to compute
  the diversity product (sum) it is enough to do $L-1$ distance
  computations.
\end{lem}
\begin{proof}
$$
|\det(\Psi_i-\Psi_j)|=|\det(I-\Psi_i^{-1}\Psi_j)| =|\det(I-B)|,
$$
where $B\in \V$ is an element in $\V$ equivalent to
$\Psi_i^{-1}\Psi_j$. This shows the result for the diversity
product. If one is concerned with the diversity sum then the same
argument still holds if the absolute value of the determinant $|
\;\det(\cdot)\; |$ is replaced by the Frobenius norm ${\| \;.\;
  \|}_F$.
\end{proof}

Based on this lemma we are interested in finite constellations
inside $G$ whose elements have a weak group structure and are all
normal. The following theorem provides a complete
characterization of all these constellations:
\begin{thm}                                  \label{mainSec3}
  Let $\V \subset G$ be a finite positive normal constellation (including identity element) with $L \geq
  3$ elements. If $\V$ has a weak group structure then
  $\V$ takes one of the following forms:
\begin{itemize}
\item $ \{I, A, A^2, \cdots, A^{L-1}\} $
\item $ \{I, AB, A^2B^2, \cdots, A^{L-1}B^{L-1}\} $
\end{itemize}

where $A=g_i^{p_i}$, $B=g_j^{p_j}$ for some $i \neq j$.
\end{thm}

The proof of Theorem~\ref{mainSec3} is rather involved. In order
to make it more  understandable we will divide it in several
definitions and lemmas.

\begin{de}
  For any element $\Psi \in G$, we define the length of
  $\Psi=\prod_{i=1}^k {a_i}^{p_i}$ to be
  $$
  \length(\Psi)=\sum_{i=1}^k p_i.
  $$
\end{de}

It is a routine to check that the definition is well-defined and
doesn't depend on the representation of the element.
For the identity element one will have $\length(I)=0$. One
immediate consequence from this definition is that if $A \sim B$, one will have
$|\length(A)|=|\length(B)|$. The following lemma claims that any
freely generated positive weak group constellation
``approximately'' takes cyclic form.

\begin{lem}
  Let $\V=\{\Psi_0=I, \Psi_1, \Psi_2, \cdots, \Psi_{L-1}\}\subset
  G$ be a positive constellation of the freely generated set
  $G\subset U(M)$. Suppose $\length(\Psi_i) \leq \length(\Psi_j)$
  for $i < j$. If $\V$ is a weak group constellation, then
  $$
  \Psi_i \in [\Psi_1]^i
  $$
  where $[\Psi_1]^i=\{a_1a_2\cdots a_i|a_1, a_2, \cdots, a_i
  \in [\Psi_1] \}$.
\end{lem}
\begin{proof}
  We first show that $\length(\Psi_i) < \length(\Psi_j)$ for $i <
  j$: Indeed, if $\length(\Psi_i) = \length(\Psi_j)$, then
  $\length(\Psi_i^{-1}
  \Psi_j)=\length(\Psi_j)-\length(\Psi_i)=0$. That means
  $\Psi_i^{-1} \Psi_j \sim I$, equivalently one will have
  $\Psi_i^{-1} \Psi_j =I$, i.e. $\Psi_i=\Psi_j$. That contradict
  the fact that $\Psi_i$ and $\Psi_j$ are distinct.

  Consider $\Psi_1^{-1} \Psi_2$. Since $0 <
  \length(\Psi_1^{-1}\Psi_2)=\length(\Psi_2)-\length(\Psi_1) <
  \length(\Psi_2)$, therefore $\Psi_1^{-1} \Psi_2=\bar{\Psi}_1$
  where $\bar{\Psi}_1 \sim \Psi_1$. So $\Psi_2=\Psi_1
  \bar{\Psi}_1 \in [\Psi_1]^2$. Proceed by induction, one
  can show $\Psi_k^{-1} \Psi_{k+1}= \bar{\Psi}_2$ where
  $\bar{\Psi}_2 \sim \Psi_1$. So $\Psi_{k+1}=\Psi_k \bar{\Psi}_2
  \in [\Psi_1]^{k+1}$ by induction.
\end{proof}

\begin{rem}
  An immediate observation is that
  $$
  \length(\Psi_i)=i * \length(\Psi_1).
  $$
\end{rem}

Take two positive normal elements in $G$ with their reduced
forms:
$$
\Psi_1=a_1^{p_1} a_2^{p_2} \cdots a_m^{p_m} \qquad
\Psi_2=b_1^{q_1} b_2^{q_2} \cdots b_n^{q_n}.
$$
We define the shift operator $S_k$ on the reduced form of a positive
normal element $\Psi$ by induction: $S_1(\Psi)=S_1(a_1^{p_1}
a_2^{p_2} \cdots a_m^{p_m})=a_2^{p_2} \cdots a_m^{p_m} a_1^{p_1}$
and $S_{k+1}=S_k \circ S_1$. We assume that $S_0(\Psi)=\Psi$,
then apparently for a fixed element $\Psi$ shift operator is
periodic. We have the following lemma.
\begin{lem} \label{shift-version}
  $\Psi_1 \sim \Psi_2$ if and only if $\Psi_1=S_k(\Psi_2)$ for
  some $k$.
\end{lem}
\begin{proof}
  The sufficiency part of this lemma is straightforward.  So we
  have to prove the necessity part. Since $\Psi_1 \sim \Psi_2$,
  according to the definition of equivalence there exists $c$
  such that $c\Psi_1 c^{-1}=\Psi_2$ or $c\Psi_1
  c^{-1}=\Psi_2^{-1}$.  However since $\length(c\Psi_1
  c^{-1})=\length(\Psi_2) > 0$ and $\length(\Psi_2^{-1}) < 0$, the
  second case won't happen. The only possibility is $c\Psi_1
  c^{-1}=\Psi_2$. We assume that $c$ is generated by only one
  generator and further assume $c=c_1^{l_1}$ with $l_1 > 0$, then
  we will have
  $$
  c_1^{l_1} a_1^{p_1} a_2^{p_2} \cdots a_m^{p_m} c_1^{-l_1}=
  b_1^{q_1} b_2^{q_2} \cdots b_n^{q_n}.
  $$
  So $c_1=a_m$ and $l_1 \leq p_m$ follows, otherwise the left
  hand side of the equation above will have negative power, while
  the right hand side only has positive power. This will
  contradict the uniqueness of the representation of the same
  element. In fact $l_1=p_m$, since otherwise $\Psi_2=c_1^{l_1}
  a_1^{p_1} a_2^{p_2} \cdots c_1^{p_m-l_1}$. This will contradict
  the fact that $\Psi_2$ is a normal element. So with
  $$
  a_m^{p_m} a_1^{p_1} \cdots a_{m-1}^{p_{m-1}}=b_1^{q_1}
  b_2^{q_2} \cdots b_n^{q_n},
  $$
  one can check $m=n$ and $\Psi_2=S_{m-1}(\Psi_1)$.

  Proceed by induction, suppose $c$ has the reduced form
  $c=c_1^{l_1} c_2^{l_2} \cdots c_{k+1}^{l_{k+1}}$, then the
  following equation follows:
  $$
  c_1^{l_1} c_2^{l_2} \cdots c_{k+1}^{l_{k+1}} a_1^{p_1}
  a_2^{p_2} \cdots a_m^{p_m} c_{k+1}^{-l_{k+1}} \cdots c_2^{-l_2}
  c_1^{-l_1}= b_1^{q_1} b_2^{q_2} \cdots b_n^{q_n}.
  $$
  Without loss of generality, we assume $l_{k+1} > 0$ and
  apply the same argument as in the one generator case. One proves
  $a_m=c_{k+1}$ and $l_{k+1}=p_m$. Therefore we reach the
  following equation:
  $$
  c_1^{l_1} c_2^{l_2} \cdots c_{k}^{l_{k}} S_{m-1}(\Psi_1)
  c_{k}^{-l_{k}} \cdots c_2^{-l_2} c_1^{-l_1}= b_1^{q_1}
  b_2^{q_2} \cdots b_n^{q_n}.
  $$
  By induction, $\Psi_2=S_{k_1} \circ
  S_{m-1}(\Psi_1)=S_{k_1+m-1}(\Psi_1)$ for some $k_1$.
\end{proof}

\begin{proof}[Proof of Theorem~\ref{mainSec3}]
  Pick any two distinct elements $ \Psi_i, \Psi_j \in \V$ having
  $\length(\Psi_i) < \length(\Psi_j)$. We claim that if $\Psi_i=a_1 a_2
  \cdots a_m$, then either there exists $1 \leq k \leq m-1$ such
  that $\Psi_j=a_1 a_2 \cdots a_k b_1 b_2 \cdots b_l a_{k+1}
  \cdots a_m$, or $\Psi_j=b_1 b_2 \cdots b_l a_1 a_2 \cdots a_m$
  or $\Psi_j=a_1 a_2 \cdots a_m b_1 b_2 \cdots b_l $ for some $l
  > 0$.

  Suppose that the claim is not true, then for $\Psi_j=c_1 c_2,
  \cdots c_p$, there exist $k_1, k_2$ such that $0 \leq k_1 \leq
  m$, $1 \leq k_2 \leq m+1$ and $k_1 < k_2-1$ and $\Psi_j$ will
  take the following form:
  $$
  \Psi_j=a_1 a_2 \cdots a_{k_1} b_1 b_2 \cdots b_l a_{k_2}
  \cdots a_m,
  $$
  where $b_1 \neq a_{k_1+1}$ and $b_l \neq a_{k_2-1}$. (For
  the special case $k_1=0$, we assume $c_1 \neq a_1$. For the
  special case $k_2=m+1$, we assume $c_p \neq a_m$.) Then
  $\Psi_i^{-1} \Psi_j$ would be equivalent to $a_{k_2-1}^{-1}
  \cdots a_{k_1+1}^{-1} b_1 b_2 \cdots b_l$, which in any case
  won't be equivalent to any positive element $\Psi_k=d_1 d_2
  \cdots d_q$ or $I$. That contradicts the fact that $\V$ is
  equipped with a weak group structure.

  As explained above we can further assume that
  $$
  \length(I) < \length(\Psi_1) < \cdots < \length(\Psi_{L-1}).
  $$

  If $\Psi_1$ is generated by only one generator, i.e.
  $\Psi_1=g_i^{p_i}$ for some $i$. Since $\Psi_2$ is a normal
  element, according to the claim, either $\Psi_2=\Psi_1
  \tilde{\Psi}_2$ or $\Psi_2=\tilde{\Psi}_2 \Psi_1$ for some
  $\tilde{\Psi}_2$. In either case $\tilde{\Psi}_2$ will be
  equivalent to $\Psi_1$, while Lemma~\ref{shift-version} will
  guarantee $\tilde{\Psi}_2=\Phi_1$. Therefore we will have
  $\Psi_2=g_i^{2p_i}$. Proceed by induction, it can be checked
  that $\Psi_l=g_i^{lp_i}$ for every $l$. So the constellation
  will take the first form in the theorem.

  If $\Phi_1$ is generated by two generators, i.e.
  $\Psi_1=g_i^{p_i} g_j^{p_j}$ for some $i, j$. According to the
  claim, we will have $\Psi_2=\Psi_1 \tilde{\Psi}_2$ or
  $\Psi_2=\tilde{\Psi}_2 \Psi_1 $ or $\Psi_2=g_i^{p_i}
  \tilde{\Psi}_2 g_j^{p_j}$. Because $\tilde{\Psi}_2$ is
  equivalent to $\Psi_1$, $\tilde{\Psi}_2$ is a shifted version
  of $\Psi_1$. Exhausting all the possibilities, the first two
  cases would make $\Psi_2$ a non-normal element, so the only
  possibility is the third case. Consider two shifted version of
  $\Psi_1$: $S_0(\Psi_1)=g_i^{p_i} g_j^{p_j} $ and
  $S_1(\Psi_1)=g_j^{p_j} g_i^{p_i}$. Only $S_0(\Psi_1)$ will
  satisfy the condition that $\Psi_2$ is a normal element. So the
  analysis above shows that
  $$
  \Psi_2=g_i^{p_i} \Psi_1 g_j^{p_j}= g_i^{2 p_i} g_j^{2 p_j}.
  $$
  By induction it can shown that
  $$
  \Psi_{k+1}=g_i^{p_i} \Psi_k g_j^{p_j}= g_i^{(k+1) p_i}
  g_j^{(k+1) p_j}.
  $$
  So in this case, the constellation will take the second form
  in the theorem.

  However the constellation doesn't exist if $\Psi_1$ is
  generated by more than $3$ elements. Indeed suppose with the
  reduced form $\Psi_1=a_1^{p_1} a_2^{p_2} \cdots a_m^{p_m}$ with
  $m \geq 3$, then $\Psi_2$ will take one of the following form:
  $\tilde{\Psi}_2 a_1^{p_1} a_2^{p_2} \cdots a_m^{p_m}, a_1^{p_1}
  \tilde{\Psi}_2 a_2^{p_2} \cdots a_m^{p_m}, \cdots, a_1^{p_1}
  a_2^{p_2} \cdots a_m^{p_m} \tilde{\Psi}_2$ with
  $\tilde{\Psi}_2$ being a shifted version of $\Psi_1$. But
  $\Psi_2$ wouldn't be a normal element for any of the above
  form, so there doesn't exist weak group constellation for this
  case.

\end{proof}

A weak group constellation is very group like, while it is not
exactly a group. It does keep the advantage of a group
constellation: for example, for any weak group constellation $\V$
taking the second form in the theorem, only $L-1$ computations
$|\det(I-A^kB^k)|$ for $k=1,2,\cdots,L-1$ are needed to calculate
the diversity product. It also overcome the disadvantage of group
codes: one can freely choose the generators, while in group
structures, the generators have to satisfy certain
relations to be a group. Last but not least it turns out that the
restriction to code elements in normal form is very advantageous
during sphere decoding. In the next section we will mainly
use the second weak group structure as described in
Theorem~\ref{mainSec3}. Before we describe these search
procedures we would like to illustrate some alternative methods.

It is possible to increase the number of generators to obtain new
structures. For instance, $\V = \{A^kB^lC^m|A, B, C \in U(M), k=0,
\cdots, p, l=0, \cdots, q, m=0, \cdots, r \}$.

For a unitary constellation $\V=\{\Phi_i|i=1, \cdots, L\}$, we
call $\V_s=\{U\Phi_iV|i=1,\cdots,L\}$ shifted version of $\V$. It
will be straightforward to prove that $\V_s$ has the same
complexity as $\V$ when one calculates the diversity.
$\{A^kCB^k|A, B, C \in U(M), k=0, \cdots, L-1\}$ is a shifted
copies of the second weak group structure in Theorem~\ref{mainSec3}. To see this,
note that $A^kCB^k=A^kCB^kC^{-1}C=A^k(CBC^{-1})^kC$. It can
checked that $A^kB^{L+1-k}=A^k{(B^{-1})}^kB^{L+1}$, therefore
$\{A^kB^{L+1-k}|A, B \in U(M), k=1, \cdots, L\}$ is also a
shifted version of the second form weak group structure.

Also we can consider the ``combination'' or the ``product'' of
two structures. For example, $\{I, A, AB, ABA, ABAB, ABABA,
\cdots\}$ is the union of $\{(AB)^k|k=0,\cdots\}$ and its shifted
version $\{(AB)^kA|k=0,\cdots\}$. Another example is the product
case: let $\V_1=\{I, C, C^2, C^3, \cdots \}$ and $\V_2=\{I, A,
AB, ABA, \cdots\}$ and consider the Cartesian product
constellation
$$
\V = \V_1 \times \V_2=\{AB|A \in \V_1, B\in \V_2\}.
$$

One may wonder how restrictive the proposed structures are. We all
know a compact Lie group can be generated by any open
neighborhood of any element in the Lie group. So with the above
structure, even if one chooses the generators locally, the
elements in the constellation could be spreading out on the whole
manifold. Somehow this indicates that the proposed structure
won't be too restrictive.

%%%%%%%%%%%%%%%%%%%%%%%%%%%%%%%%%%%%%%%%%%%%%%%%%%%%%%%%%%%%%%%%
\Section{Geometrical Design of Unitary Constellations with Good
  Diversity}     \label{Sec-geometrical}

For low dimensional constellations, one may further specify the
generators in the proposed structure. Observe that for the second
form weak group constellation, one can always assume $A$ is
diagonal. In the sequel, we further assume that $B$ is real
orthogonal, i.e.  based on the weak group structure we consider
the following $2$ dimensional constellation:
\begin{equation}  \label{Closed-Specification}
\V=\{  A^kB^k| A=\left( \begin{array}{cc}
                e^{ix}&0\\
                0&e^{iy}
          \end{array} \right), B=\left( \begin{array}{cc}
                \cos z& \sin z\\
                -\sin z& \cos z
          \end{array} \right), k=0,1,\cdots,L-1\}.
\end{equation}

There are several ways to design constellations with good diversity
from this specific structure. A natural idea is to do
Brute Force search using fine step size. Another approach is to
design the constellation with the help of geometrical intuition.
Note that a $2 \times 2$ complex matrix can be viewed as a vector
in $\mathbb{C}^4$. In this context $A$ and $B$ can be viewed as
``rotation'' transforms (induced by regular matrix multiplication)
acting on $\mathbb{C}^4$. A constellation of
form\eqr{Closed-Specification} can be viewed as a set of rotated
vectors under the transforms $A^kB^k$, $k=0, 1, \cdots, L-1$.
Intuition says that good constellations can be found if the
rotation angle is symmetrical. Based on the idea above we assume that
$x, y, z$ to be the multiples of $2\pi/L$, we found a lot of good
codes resulted from this geometrical symmetry (see tables in
Section~\ref{Sec-numerical}).

\begin{figure}[ht]
\centerline{\psfig{figure=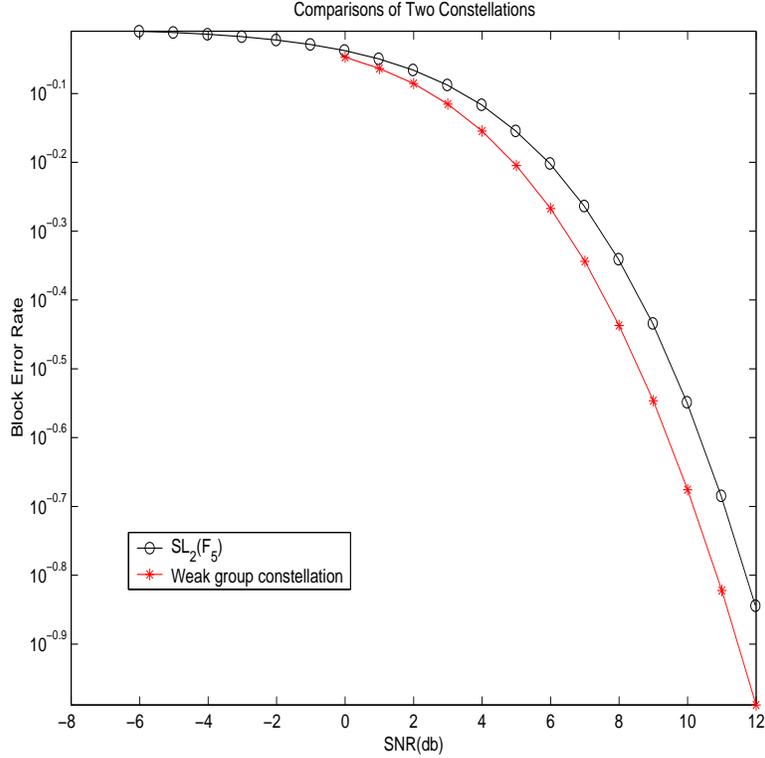,width=4in,height=4in}}
\caption{2 dimensional weak group constellations and group
constellation} \label{Closed-Group}
\end{figure}

2 dimensional constellation design has been studied
in~\cite{li02}. In this paper Liang proposed very interesting parametric
codes and many codes with excellent diversity are found. The
codes shown in~\cite{li02} can be achieved by our design as
well. In fact, most of Liang's codes belong to a special form of our
parameterization\eqr{Closed-Specification}. To our best
knowledge, most of our codes shown on the web site~\cite{ha03u2}
are the best codes ever found or never found before.

\begin{exmp}
A very interesting code with $120$ elements is found using this
approach:
$$
\V=\{ A^kB^k  |A=\left( \begin{array}{cc}
                e^{\pi/30 i}&0\\
                0&e^{11\pi/30 i}
          \end{array} \right), B=\left( \begin{array}{cc}
                \cos \pi/4& \sin \pi/4\\
                -\sin \pi/4& \cos \pi/4
\end{array} \right), k=0,1,\cdots,119\}.
$$

It can be checked that $\prod {\V}=\sum {\V} =
\frac{1}{2}\sqrt{\frac{(3-\sqrt{5})}{2}}$, i.e. the diversity
product and the diversity sum are identical to the ones of
the $SL_2(\mathbb{F}_5)$-constellation. We
simulated the performance of this code and compared it with the
performance of the $SL_2(\mathbb{F}_5)$-constellation. To our big
surprise our new code performed considerably better than the
$SL_2(\mathbb{F}_5)$-constellation. The constellation $\V$ with
sphere decoding outperformed the
$SL_2(\mathbb{F}_5)$-constellation by about $1db$ up to about
$20$db (see Figure~\ref{Closed-Group}). As the SNR goes higher,
the two curves are coming closer though.

In order to understand the difference in the performance of the
two seemingly similar constellations we investigated the
diversity product (DP) and diversity sum (DS) {\em distance
  spectrum} for each of them. As we explained before, for a
unitary constellation with $L$ elements, $L(L-1)/2$ distance
calculations may produce distances with multiplicities. For
example consider $\V$ as above, $360$ out of $7140$ pairs of
elements have distance $0.3090$ (see DP distance spectrum in
Table~\ref{DP-DS}). So one can explain the behavior difference of
the two codes using their distance spectrum. The following
Table~\ref{DP-DS} shows that the DP and DS distance spectrum of
our weak group constellation.

\begin{ta}                \label{DP-DS}
$$
\begin{array}{cc}
\begin{tabular}{c}
  Weak group constellation \\
  DP distance spectrum \\
\end{tabular} & \begin{tabular}{c}
  Weak group constellation \\
  DS distance spectrum \\
\end{tabular}\\
   \begin{tabular}{|c|c|}
     \hline
     distance & distribution \\
     \hline
     0.3090 & 360 \\
     \hline
     0.3136 & 480 \\
     \hline
     0.3895 & 480 \\
     \hline
     0.3931 & 1440 \\
     \hline
     0.4402 & 240 \\
     \hline
     0.5000 & 120 \\
     \hline
     0.5878 & 120 \\
     \hline
     0.6360 & 1440 \\
     \hline
     0.6787 & 480 \\
     \hline
     0.7071 & 600\\
     \hline
     0.8090 & 360\\
     \hline
     0.8430 & 480 \\
     \hline
     0.8660 & 120\\
     \hline
     0.8979 & 240 \\
     \hline
     0.9511 & 120 \\
     \hline
     1 & 60 \\
     \hline
   \end{tabular}& \begin{tabular}{|c|c|}
     \hline
     distance & distribution\\
     \hline
     0.3090 & 120 \\
     \hline
     0.4402 & 240 \\
     \hline
     0.5000 & 120 \\
     \hline
     0.5023 & 480 \\
     \hline
     0.5457 & 240 \\
     \hline
     0.5878 & 120 \\
     \hline
     0.6367 & 480 \\
     \hline
     0.6502 & 240 \\
     \hline
     0.7071 & 3000 \\
     \hline
     0.7598 & 240 \\
     \hline
     0.7711 & 240 \\
     \hline
     0.8090 & 120 \\
     \hline
     0.8380 & 240 \\
     \hline
     0.8647 & 480 \\
     \hline
     0.8660 & 120 \\
     \hline
     0.8979 & 240 \\
     \hline
     0.9511 & 120 \\
     \hline
     1 & 60 \\
     \hline
   \end{tabular}
\end{array}
$$
\end{ta}
One can check that the DP distance spectrum of the
$SL_2(\mathbb{F}_5)$-constellation is identical to the DS
distance spectrum. The following Table~\ref{DP-DP} shows that the DS distance
spectrum for the $SL_2(\mathbb{F}_5)$-constellation has denser
small distance distribution compared to DS spectrum of our
constellation and this explains the considerable worse
performance of this constellation in our simulations.
\begin{ta}          \label{DP-DP} \
\begin{center}
\begin{tabular}{c}
  $SL_2(\mathbb{F}_5)$-constellation \\
  DP (DS) distance spectrum \\
\end{tabular} \\
\begin{tabular}{|c|c|}
     \hline
     distance & distribution \\
     \hline
     0.3090 & 720 \\
     \hline
     0.5000 & 1200 \\
     \hline
     0.5878 & 720 \\
     \hline
     0.7071 & 1800 \\
     \hline
     0.8090 & 720 \\
     \hline
     0.8660 & 1200 \\
     \hline
     0.9511 & 720 \\
     \hline
     1 & 60 \\
     \hline
\end{tabular}
\end{center}
\end{ta}
\end{exmp}

Although we have concentrated so far in the design of
2-dimensional constellations there is actually no restriction
with our approach. The similar ``rotation'' idea can be applied
to other low dimensional constellation design. For instance, we
can make further specifications to $3$ dimensional weak group
constellations:
$$
\V=\{ A^kB^k |A=\left( \begin{array}{ccc}
                \cos x & \sin x & 0\\
                -\sin x& \cos x & 0\\
                0      &   0    & e^{iy}
          \end{array} \right), B=\left( \begin{array}{ccc}
                e^{iz} & 0 & 0\\
                0      & \cos w & \sin w \\
                0      & -\sin w & \cos w
          \end{array} \right), k=0,1,\cdots,L-1\}.
$$
where $x,y,z,w$ is assumed to take the multiple of $2\pi/L$.
Apparently algebraic design based on geometrical symmetry can be
applied to any other structure as well. For instance consider the
following specified structures:
$$
\V=\{  A^kB^l| A=\left( \begin{array}{cc}
                e^{ix}&0\\
                0&e^{iy}
          \end{array} \right), B=\left( \begin{array}{cc}
                \cos z& \sin z\\
                -\sin z& \cos z
          \end{array} \right), k=0,1,\cdots,p-1, l=0,1,\cdots,q-1\}.
$$
where we can take $x,y$ to be multiple of $2\pi/p$ and $z$ to be
multiple of $2\pi/q$.  We refer to~\cite{ha03u2} for the designed
low dimensional constellations from these approaches.

%%%%%%%%%%%%%%%%%%%%%%%%%%%%%%%%%%%%%%%%%%%%%%%%%%%%%%%
\Section{Numerical Design of Unitary Constellation with Good
Diversity}  \label{Sec-numerical}

In order to numerically design constellations, it will be necessary to have a
good parameterization for the set of unitary constellations having
size $L$, operating with $M$ transmit antennas. In this section we
show how one can use the theory of complex Stiefel manifolds and
the classical Cayley transform to obtain such a parameterization.

\subsection{The complex Stiefel manifold}

\begin{de}
  The subset of $T\times M$ complex matrices
  $$
  \S_{T,M}:=\left\{ \Phi\in\C^{T\times M}\mid \Phi^* \Phi =
    I_M\right\}
  $$
  is called the {\em complex Stiefel manifold}.
\end{de}

{}From an abstract point of view a constellation $\V:=\{
\Phi_1,\ldots, \Phi_L\}$ having size $L$, block length $T$ and
operating with $M$ antennas can be viewed as a point in the
complex manifold
$$
\mathcal{M}:=\left(\S_{T,M}\right)^L=
\underbrace{\S_{T,M}\times\cdots \times\S_{T,M}}_{\mbox{$L$
    copies}}.
$$
The search for good constellations $\V$ requires hence the
search for points in $\mathcal{M}$ whose diversity is excellent
in some interval $[\rho_1,\rho_2]$.

Stiefel manifolds have been intensely studied in the mathematics
literature since their introduction by Eduard Stiefel some 50
years ago. A classical paper on complex Stiefel manifolds
is~\cite{at60}, a paper with a point of view toward numerical
algorithms is~\cite{ed99}. The major properties are summarized by
the following theorem:
\begin{thm}                                     \label{Stiefel}
  $\S_{T,M}$ is a smooth, real and compact sub-manifold of
  $\C^{MT}=\R^{2MT}$ of real dimension $2TM-M^2$.
\end{thm}

Some of the stated properties will follow from our further
development. The following two examples give some special cases.

\begin{exmp}
  $$
  \S_{T,1}=\left\{ x\in\C^T\mid ||x||=\sqrt{\sum_{i=1}^M
      x_i\bar{x}_i}=1 \right\}\subset \R^{2T}
  $$
  is isomorphic to the $2T-1$ dimensional unit sphere
  $S^{2T-1}$.
\end{exmp}

\begin{exmp}
  When $T=M$ then $\S_{T,M}=U(M)$, the group of $M\times M$
  unitary matrices. It is well known that the Lie algebra of
  $U(M)$, i.e. the tangent space at the identity element,
  consists of all $M\times M$ skew-Hermitian matrices. This
  linear vector space has real dimension $M^2$, in particular the
  dimension of $U(M)$ is $M^2$ as well.
\end{exmp}

A direct consequence of Theorem~\ref{Stiefel} is:

\begin{co}
  The manifold $\mathcal{M}$ which parameterizes the set of all
  constellations $\V$ having size $L$, block length $T$ and
  operating with $M$ antennas forms a a real compact manifold of
  dimension $2LTM-LM^2$.
\end{co}

As this corollary makes it clear a full search over the total
parameter space is only possible for very moderate sizes of
$M,L,T$. It is also required to have a good parameterization of
the complex Stiefel manifold $\S_{T,M}$ and we will go after this
task next.

The unitary group is closely related to the complex Stiefel
manifold and the problem of parameterization ultimately boils
down to the parameterization of unitary matrices. For this assume
that $\Phi$ is a $T\times M$ matrix representing an element of
the complex Stiefel manifold $\S_{T,M}$. Using Gramm-Schmidt one
constructs a $T\times (T-M)$ matrix $V$ such that the $T\times T$
matrix $\left[ \Phi\mid V\right]$ is unitary. Define two $T\times
T$ unitary matrices $\left[ \Phi_1\mid V_1\right]$ and $\left[
  \Phi_2\mid V_2\right]$ to be equivalent whenever
$\Phi_1=\Phi_2$. A direct calculation shows that two matrices are
equivalent if and only if there is $(T-M)\times (T-M)$ matrix $Q$
such that:
\begin{equation}                                       \label{Q-matrix}
\left[ \Phi_2\mid V_2\right]=\left[ \Phi_1\mid V_1\right]\vier{I}{0}{0}{Q}.
\end{equation}
Identifying the set of matrices $Q$ appearing in\eqr{Q-matrix}
with the unitary group $U(T-M)$ we get the result:
\begin{lem}                                   \label{Lem-par}
  The complex Stiefel manifold $\S_{T,M}$ is isomorphic to the
  quotient group
  $$
  U(T)/U(T-M).
  $$
\end{lem}
This lemma let us verify the dimension formula for $\S_{T,M}$
stated in Theorem~\ref{Stiefel}:
$$
\dim \S_{T,M}=\dim U(T)-\dim U(T-M)=T^2-(T-M)^2=2TM-M^2.
$$

The  section makes it clear that a good parameterization of the
set of constellations $\V$ requires a good parameterization of the
manifold $\mathcal{M}$ and this in turn requires a good
parameterization of the unitary group $U(M)$.

Once one has a nice parameterization of the unitary group $U(M)$
then Lemma~\ref{Lem-par} provides a way to parameterize the
Stiefel manifold $\S_{T,M}$ as well. Parameterizing $U(T)$ modulo
$U(T-M)$ is however an `over parameterization'. Edelman, Arias and
Smith~\cite{ed99} explained a way on how to describe a local
neighborhood of a (real) Stiefel manifold $\S_{T,M}$. The method
can equally well be applied in the complex case. We do not pursue
this parameterization in this paper and leave this for future
work.

In the remainder of this paper we will concentrate on
constellations having the special form\eqr{specialform}. From a
numerical point of view we require for this a good
parameterization of the unitary group and the next subsection
provides an elegant way to do this.

%%%%%%%%%%%%%%%%%%%%%%%%%%%%%%%%%%%%%%%%%%%%%%%%%%%%%%%%
\subsection{Cayley transformation}

There are several ways to represent a unitary matrix in a very
explicit way. One elegant way makes use of the classical Cayley
transformation. In order that the paper is self contained we
provide a short summary. More details are given
in~\cite[Section~22]{pr94} and~\cite{ha02a}.

\begin{de}
  For a complex $M \times M$ matrix $Y$ which has no eigenvalues at $-1$,
  the Cayley transform of $Y$ is defined to be
  $$
  Y^c= (I + Y)^{-1}(I-Y),
  $$
  where $I$ is the $M \times M$ identity matrix.
\end{de}
Note that $(I+Y)$ is nonsingular whenever $Y$ has no eigenvalue
at -1. One immediately verifies that $(Y^c)^c=Y$. This is in
analogy to the fact that the linear fractional transformation
$f(z)=\frac{1-z}{1+z}$ has the property that $f(f(z))=z$. Recall
that a matrix $M$ is skew-Hermitian whenever $A^*=-A$. The set of
$M\times M$ skew-Hermitian matrices forms a linear subspace of
$\C^{M\times M}\cong \R^{2M^2}$ having real dimension $M^2$. This
is the Lie algebra of the unitary group $U(M)$. The main property
of the Cayley transformation is summarized in the following
theorem. (See e.g.~\cite{ha02a,pr94}).
\begin{thm}
  When $A$ is a skew-Hermitian matrix then $(I+A)$ is nonsingular
  and the Cayley transform $V:=A^c$ is a unitary matrix. Vice versa
  when $V$ is a unitary matrix which   has no eigenvalues at $-1$
  then the Cayley transform $V^c$ is skew-Hermitian.
\end{thm}

This theorem allows one to parameterize the open set of $U(M)$
consisting of all unitary matrices whose eigenvalues do not
include $-1$ through the linear vector space of skew-Hermitian
matrices. The Cayley transformation is very important for the
numerical design of constellations because it makes the local
topology of $U(M)$ clear. One can see that most optimization
method require us to consider the neighborhood of one element in
$U(M)$.

%%%%%%%%%%%%%%%%%%%%%%%%%%%%%%%%%%%%%%%%%%%%%%%%%%%%%%%%%%%%%%%%%%%%%
\subsection{Simulated Annealing (SA) Algorithm}

In our numerical experiments we have considered several methods.
Because there are a large number of target functions the best
known optimization algorithms such as Newton's
Methods~\cite{no99,ed99} and the Conjugate Gradient
Method~\cite{no99,ed99} are difficult to implement. Surprisingly
the {\em Simulated Annealing Algorithm} turned out to be very
practical for this problem.

Simulated Annealing (SA) is a method which mimics the process of
melted metal getting cooled off. In the annealing process of the
melted metal, first the metal is heated to melt, then the
temperature is getting down gradually. The metal will get to a
minimized energy state if the temperature is lowering slow enough.
For more details about this algorithm, we refer
to~\cite{aa89,la87b,ot89}.

In fact, we would rather call it a general method instead of a
concrete algorithm. Generally speaking, for a given optimization
problem we always take an initial solution in some certain way,
then consider a second solution in the ``neighborhood'' of this
solution. We will accept the solution according to some predefined
criterion which might involve a probability threshold.

Combining with good algebraic structure and Cayley transform,
which is a good representation of any dimensional unitary matrix,
one can see that numerical method can be applied to any
dimensional and any size constellation design. Our implementation
of the algorithm can be summarized in the following way, one can
find simple sample program on our web site~\cite{ha03u2}.

\begin{enumerate}

\item Choose a proposed algebraic structure for the constellation.

\item Generate initial generators of the whole constellation.
      One can either take an existing constellation as the start
      point or just take the initial point randomly.

\item Generate randomly a new constellation using Cayley transform in
  the neighborhood of the old constellation where the selection is done using a Gaussian
  distribution with decreasing variances as the algorithm
  progresses.

\item Calculate the diversity function (product, sum) of the
  newly constructed constellation.

\item If the new constellation has better diversity function
  (product, sum), then accept the new constellation. If not,
  reject the new constellation and keep the old constellation (or
  accept it according to Metropolis's criterion~\cite{me53}).

\item Check the stopping criterion, if satisfied, then stop,
  otherwise go to $2$ and continue the iteration.

\end{enumerate}

\begin{exmp}
As we mentioned before, one can either choose an existing
constellation as the starting point for our numerical method or
just take the initial point randomly. In the sequel, we use the
group constellation $G_{21,4}$ in~\cite{sh01}:

$$
\V_1=\{A^kB^l| A=\left(\begin{array}{ccc}
         \eta&0&0\\
         0&\eta^4&0\\
         0&0&\eta^{16}\\
     \end{array}\right), B=\left(\begin{array}{ccc}
                                  0&1&0\\
                                  0&0&1\\
                                  \eta^7&0&0\\
                                  \end{array}\right),
                                  k=0,1,\cdots,20, l=0,1,2\}
$$

One can verify that
$$\prod \V_1=0.3851.$$
It seems that $G_{21,4}$ is already a very good constellation,
our algorithm only improves a little (see $\V_2$ below). However
one can check for most of the cases, the algorithm will improve
much compared to the original group constellation.
$$\V_2=\{A^kB^l|k=0,1,\cdots,20,l=0,1,2 \},$$
where
$$A=\left(\begin{array}{ccc}
         0.9415 + 0.3155*i&0.0573 - 0.0222*i&0.0496 +0.0882*i\\
         0.0160 - 0.0555*i&0.4005 + 0.9136*i&0.0326 - 0.0212*i\\
         0.0579 + 0.0855*i&-0.0312 - 0.0099*i&0.1384 - 0.9844*i\\
                 \end{array}\right),$$
$$ B=\left(\begin{array}{ccc}
        0.0175 + 0.0095*i&0.9997 + 0.0111*i&0.0079 + 0.0042*i\\
        0.0086 + 0.0100*i&-0.0082 + 0.0040*i&0.9999 + 0.0036*i\\
        -0.4836 + 0.8750*i&0.0004 - 0.0198*i&-0.0045 - 0.0126*i\\
        \end{array}\right).$$
One verifies that
$$
\prod \V_2=0.3874.
$$
\end{exmp}

\begin{exmp}
Different industrial applications require different level of
reliability of the communication channels. One may want to
optimize the constellation at certain Block Error Rate (BER) or
Signal Noise Ratio (SNR). It can also be shown theoretically that
numerical methods together with the proposed structure works in
the same way if one wants to optimize the diversity function at a
certain SNR. This is essential the case because for a complex
matrix $A$ and unitary matrices $U, V$ one has that
\begin{equation} \label{sing}
\delta_m(UAV)=\delta_m(A),
\end{equation}
for $m=1,2,\cdots,M$. With the constellation structures as above
we are able to reduce the dimension of the parameter space and at
the same time we have a considerable reduction in the number of
targets to be checked. Intuitively algebraically designing codes
for this purpose seems to be impossible.

The following graph shows the comparison of three constellations
with different dimensions with $2$ receiver antennas. The first
one is a $2$ dimensional constellation with $3$ elements
($R=0.7925$) and optimal diversity product $0.8660$ and optimal
diversity sum $0.8660$.  The second constellation is a $3$
dimensional constellation which has $5$ elements ($R=0.7740$)
with diversity product $0.7183$ and diversity sum $0.7454$. The
third constellation is a $4$ dimensional one consisting of $9$
elements ($R=0.7925$) with diversity product 0.5904 and diversity
sum 0.6403. Here based on the structure $A^kB^k$ we used Simulated
Annealing to optimize the diversity function at $6$db to acquire
the last two constellations.

One can see that around $5$ db, the second constellation surpasses
the first one and is getting better and better as the SNR becomes
larger. This can be easily understood since the diversity function
of the first constellation is approximately dominated by
$1/{\rho^4}$ at high SNR, while the diversity function of the
second constellation is  dominated by $1/{\rho^6}$. The same
explanation can be applied to the third constellation's
performance. One can even foresee that higher dimensional
constellations will perform even better and the BER curve will be
sharper than the lower dimensional ones. It is believable that higher
dimensional constellations will achieve much more diversity gain
compared to lower dimensional ones.

\begin{figure}[ht]
\centerline{\psfig{figure=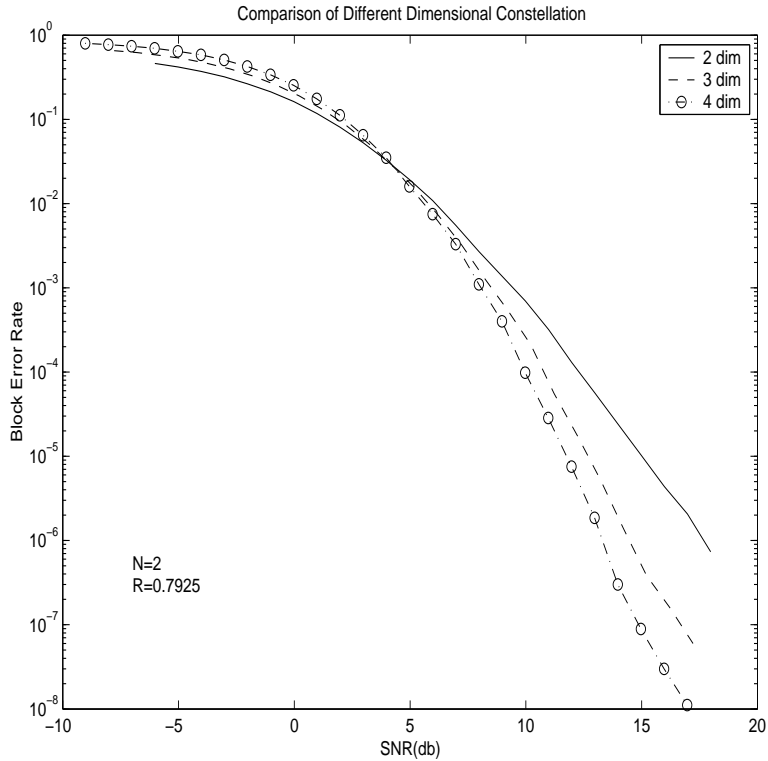,width=4in,height=4in}}
\caption{Performance of different dimensional constellations with
the same rate} \label{fig-5}
\end{figure}
\end{exmp}

Surprisingly SA works very well when it is applied to an
algebraic structure with symmetry. Like all the other numerical
methods, one has to suffer the loss of performance due to the
increasing complexity as the size and dimension go up and due to
the limited computational resources. However without any doubts,
the numerical approach is very flexible and can be used for any
dimensional and any size constellation, producing very good
diversity. So a lot of good-performing unitary constellations are
found this way, which were never found by any algebraic
method. At the end of this section we will show some $2$
dimensional constellations we found using
various methods based on the proposed structure. We skip,
however, our numerical results on the higher dimensional unitary
constellation design, since one can check them on the web
site~\cite{ha03u2}.

One very interesting fact is the numerical results for diversity
sum from $3$ dimensional structured constellation are even better
than the corresponding upper bound for $2$ dimensional
constellations. Somehow it won't be too surprising if one notices
that from $U(2)$ to $U(3)$, we have $5$ more dimensions to
manoeuvre.

In~\cite{ha03u} packing problems on compact Lie groups are
analyzed and the upper bound for the diversity sum and the
diversity product are derived. In the following figure one can
see the limiting behavior of $2$ dimensional structured
constellations compared to the upper bound. One can
check~\cite{ha03u2} for the comparisons for other dimensions.

\begin{figure}[ht]
\centerline{\psfig{figure=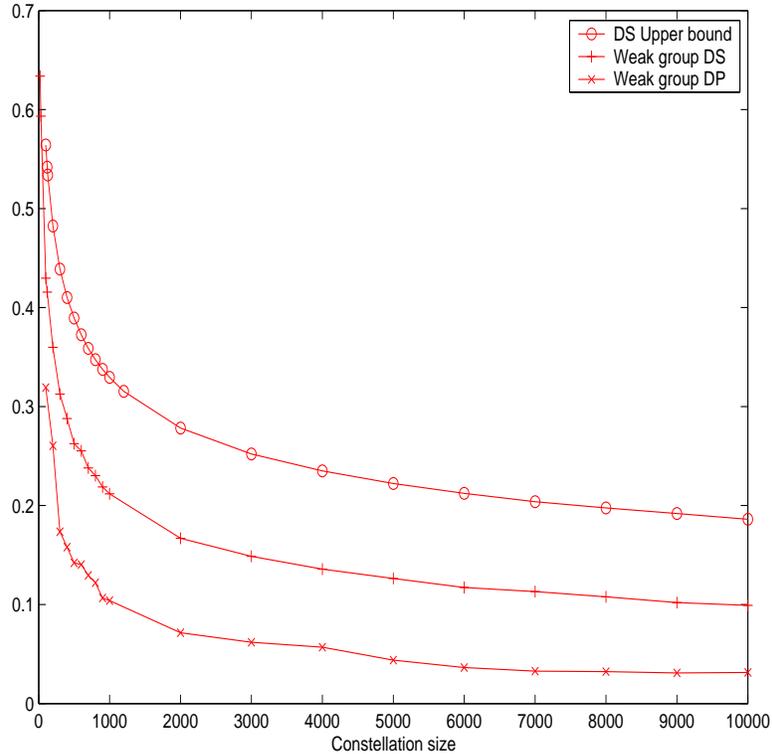,width=4in,height=4in}}
\caption{2 dimensional weak group constellations and upper bound}
\end{figure}

%%%%%%%%%%%%%%%%%%%%%%%%%%%%%%%%%%%%%%%%%%%%%%%%%%%%%%%%%%%%%%%
\subsection{Constellations with extremely large diversity}

In this subsection we list the best $2$-dimensional
constellations we found with the techniques described in
Sections~\ref{Sec-geometrical} and~\ref{Sec-numerical}. The
tabulated constellations have some of the best diversity sums and
diversity products published so far. All the constellations
searched by simulated annealing (SA) were based on the $A^kB^k$
structure. For the constellations with $L$ elements and
parameters $x, y, z$ being multiples of $2\pi/L$, they are found
by geometrical methods using the
parameterization\eqr{Closed-Specification}. For the
constellations with $L$ elements and parameters $x, y, z$ being
decimals, they are found by Brute Force with step size $0.1000$
based on the same parameterization\eqr{Closed-Specification}.

\newpage

\begin{ta}
Diversity product of $2$ dimensional constellation based on weak
group structure:

\begin{center}
\begin{tabular}{|c|c|c|}
  \hline
   \begin{tabular}{c}
   Number of\vspace*{-2mm}\\ elements\end{tabular}& Diversity Product & Codes and Comments\\
  \hline
  2 & 1 & $x=\pi, y=\pi, z=0 \; (\mbox{optimal})$\\
  \hline
  3 & $\sqrt{3}/2$  & $x=2\pi/3, y=2\pi/3,z=0 \; (\mbox{optimal})$\\
  \hline
  4 & 0.7831  &$x=0.6000, y=6.0000, z=4.4000$\\
  \hline
  5 & $\sqrt{5/8}$&$x=2\pi/5, y=8\pi/5, z=4\pi/5 \; (\mbox{optimal})$\\
  \hline
  8 & 0.7071 &$x=2.3562,y=3.9270,z=4.7124$\\
  \hline
  9 & 0.6524 & SA searched code \\
  \hline
  10& 0.6124 & $x=2\pi/5,y=8\pi/5,z=\pi/5$\\
  \hline
  16 & $\sqrt[4]{2}/2$ &$x=\pi/4,y=5\pi/4,z=13\pi/8$\\
  \hline
  17 & 0.5255 & SA searched code\\
  \hline
  18 & 0.5207 & SA searched code\\
  \hline
  19 & 0.5128 & SA searched code\\
  \hline
  20 & 0.5011 & $x=1.6500,y=3.7500,z=4.0500$\\
  \hline
  24 & 0.5000 &$x=\pi/12,y=5\pi/12,z=\pi/2$\\
  \hline
  37 & 0.4461 &$x=2\pi/37,y=6\pi/37,z=12\pi/37$\\
  \hline
  39 & 0.3984 &$x=8\pi/39,y=34\pi/39,z=36\pi/39$\\
  \hline
  40 & 0.3931 &$x=3\pi/10,y=11\pi/10,z=3\pi/4$\\
  \hline
  55 & 0.3874 &$x=2\pi/55,y=68\pi/55,z=6\pi/11$\\
  \hline
  57 & 0.3764 &$x=2\pi/57,y=40\pi/57,z=48\pi/57$\\
  \hline
  75 & 0.3535 &$x=2\pi/75,y=98\pi/75,z=96\pi/75$\\
  \hline
  85 & 0.3497 &$x=26\pi/85,y=94\pi/85,z=18\pi/17$\\
  \hline
  91 & 0.3451 &$x=2\pi/91,y=128\pi/91,z=42\pi/91$\\
  \hline
  96 & 0.3192 &$x=7\pi/16,y=29\pi/16,z=\pi/6$\\
  \hline
  105 & 0.3116 &$x=2\pi/105,y=68\pi/105,z=84\pi/105$\\
  \hline
  120 &0.3090  &$x=\pi/30,y=11\pi/30,z=\pi/4$\\
  \hline
  135 & 0.2869 &$x=2\pi/135,y=28\pi/135,z=68\pi/135$\\
  \hline
  145 & 0.2841 &$x=2\pi/145,y=64\pi/145,z=76\pi/145$\\
  \hline
  165 & 0.2783 &$x=2\pi/33,y=20\pi/33,z=2\pi/5$\\
  \hline
  203 & 0.2603 &$x=2\pi/203,y=290\pi/203,z=70\pi/203$\\
  \hline
  225 & 0.2499 &$x=82\pi/225,y=118\pi/225,z=126\pi/225$\\
  \hline
  217 & 0.2511 &$x=2\pi/217,y=250\pi/217,z=168\pi/217$\\
  \hline
  225 & 0.2499 &$x=82\pi/225,y=118\pi/225,z=126\pi/225$\\
  \hline
  240 & 0.2239 &$x=\pi/40,y=9\pi/40,z=\pi/6$\\
  \hline
  273 & 0.2152 &$x=2\pi/273,y=208\pi/273,z=142\pi/273$\\
  \hline
  295 & 0.2237 &$x=14\pi/295,y=104\pi/295,z=22\pi/59$\\
  \hline
  297 & 0.1910 &$x=242\pi/297,y=548\pi/297,z=54\pi/297$\\
  \hline
  299 & 0.1858 &$x=8\pi/299,y=220\pi/299,z=18\pi/299$\\
  \hline
  300 & 0.1736 &$x=\pi/150,y=51\pi/150,z=5\pi/6$\\
  \hline

\end{tabular}
\end{center}
\end{ta}

%\begin{tabular}{|c|c|c|}
%  \hline
%   \begin{tabular}{c}
%   number of\vspace*{-2mm}\\ elements\end{tabular}& Diversity Product & Codes and Comments\\
%  \hline
%  203 & 0.2603 &$x=2\pi/203,y=290\pi/203,z=70\pi/203$\\
%  \hline
%  225 & 0.2499 &$x=82\pi/225,y=118\pi/225,z=126\pi/225$\\
%  \hline
%  217 & 0.2511 &$x=2\pi/217,y=250\pi/217,z=168\pi/217$\\
%  \hline
%  225 & 0.2499 &$x=82\pi/225,y=118\pi/225,z=126\pi/225$\\
%  \hline
%  240 & 0.2239 &$x=\pi/40,y=9\pi/40,z=\pi/6$\\
%  \hline
%  273 & 0.2152 &$x=2\pi/273,y=208\pi/273,z=142\pi/273$\\
%  \hline
%  295 & 0.2237 &$x=14\pi/295,y=104\pi/295,z=22\pi/59$\\
%  \hline
%  297 & 0.1910 &$x=242\pi/297,y=548\pi/297,z=54\pi/297$\\
%  \hline
%  299 & 0.1858 &$x=8\pi/299,y=220\pi/299,z=18\pi/299$\\
%  \hline
%  300 & 0.1736 &$x=\pi/150,y=51\pi/150,z=5\pi/6$\\
%  \hline
%
%\end{tabular}
%\end{center}
%\end{ta}

\newpage
\begin{ta}
Diversity sum of $2$ dimensional constellation based on weak group
structure

\begin{center}
\begin{tabular}{|c|c|c|}
  \hline
   \begin{tabular}{c}
   number of\vspace*{-2mm}\\ elements\end{tabular}& Diversity Sum& Codes and Comments\\
  \hline
  2 & 1 & $x=\pi, y=\pi, z=0 \; (\mbox{optimal})$\\
  \hline
  3 & $\sqrt{3}/2$  &$x=2\pi/3, y=2\pi/3,z=0 \; (\mbox{optimal})$\\
  \hline
  5 & $\sqrt{5/8}$ &$x=2\pi/5, y=8\pi/5, z=4\pi/5 \; (\mbox{optimal})$\\
  \hline
  9 & 3/4 &$x=10\pi/9, y=4\pi/3, z=4\pi/9 \; (\mbox{optimal})$\\
  \hline
  16 & $\sqrt{2}/2$ &$x=\pi/4, y=5\pi/4, z=13\pi/8 \; (\mbox{optimal})$\\
  \hline
  18 & 0.6614 &$x=4\pi/9,y=2\pi/3,z=7\pi/9$\\
  \hline
  19 & 0.6391 & SA searched code\\
  \hline
  20 & 0.6338 & SA searched code\\
  \hline
  21 & 0.6307 & SA searched code\\
  \hline
  22 & 0.6154 & SA searched code\\
  \hline
  24 & 0.6124 & $x=\pi/6,y=\pi/4,z=5\pi/12$\\
  \hline
  28 & 0.5996 & $x=3\pi/8,y=\pi/2,z=2\pi/7$\\
  \hline
  30 & 0.5934 & $x=4\pi/15,y=\pi/3,z=7\pi/15$\\
  \hline
  31 & 0.5739 & SA searched code\\
  \hline
  32 & 0.5734 & SA searched code\\
  \hline
  39 & 0.5726 & $x=14\pi/39,y=40\pi/39,z=18\pi/39$\\
  \hline
  40 & 0.5499 & $x=3\pi/20,y=7\pi/20,z=3\pi/10$\\
  \hline
  42 & 0.5371 & $x=4\pi/7,y=13\pi/21,z=\pi/3$\\
  \hline
  45 & 0.5342 & $x=2\pi/9,y=4\pi/9,z=14\pi/15$\\
  \hline
  52 & 0.5332 & $x=\pi/13,y=2\pi/13,z=9\pi/26$\\
  \hline
  57 & 0.5053 & $x=4\pi/57,y=8\pi/57,z=40\pi/57$\\
  \hline
  60 & 0.5000 & $x=\pi/15,y=4\pi/15,z=3\pi/10$\\
  \hline
  64 & 0.4852 &$x=3\pi/16, y=53\pi/32, z=55\pi/32$\\
  \hline
  75 & 0.4850 &$x=32\pi/75,y=14\pi/75,z=2\pi/75$\\
  \hline
  76 & 0.4672 &$x=3\pi/19,y=4\pi/19,z=11\pi/38$\\
  \hline
  77 & 0.4595 &$x=52\pi/77,y=82\pi/77,z=60\pi/77$\\
  \hline
  85 & 0.4540 &$x=2\pi/17,y=8\pi/17,z=14\pi/85$\\
  \hline
  87 & 0.4460 &$x=52\pi/87,y=98\pi/87,z=82\pi/87$\\
  \hline
  95 & 0.4418 &$x=6\pi/19,y=2\pi/95,z=36\pi/95$\\
  \hline
  96 & 0.4390 &$x=39\pi/48,y=5\pi/12,z=11\pi/24$\\
  \hline
  99 & 0.4297 &$x=62\pi/99,y=192\pi/99,z=142\pi/99$\\
  \hline
  105 & 0.4295 &$x=2\pi/105,y=16\pi/105,z=28\pi/105$\\
  \hline
  106 & 0.4161 &$x=2\pi/53,y=13\pi/53,z=12\pi/53$\\
  \hline
  120 & 0.4156 &$x=\pi/10,y=\pi/6,z=5\pi/4$\\
  \hline
  123 & 0.4077 &$x=188\pi/123,y=38\pi/123,z=182\pi/123$\\
  \hline
\end{tabular}

\begin{tabular}{|c|c|c|}
  \hline
   \begin{tabular}{c}
   number of\vspace*{-2mm}\\ elements\end{tabular}& Diversity Sum& Codes and Comments\\
  \hline
  130 & 0.4071 &$x=26\pi/65,y=5\pi/13,z=2\pi/13$\\
  \hline
  133 & 0.3971 &$x=2\pi/133,y=212\pi/133,z=206\pi/133$\\
  \hline
  138 & 0.3963 &$x=16\pi/69,y=19\pi/69,z=4\pi/69$\\
  \hline
  145 & 0.3949 &$x=138\pi/145,y=22\pi/145,z=40\pi/29$\\
  \hline
  148 & 0.3840 &$x=5\pi/74,y=13\pi/37,z=2\pi/37$\\
  \hline
  150 & 0.3758 &$x=\pi/15,y=8\pi/75,z=19\pi/75$\\
  \hline
  155 & 0.3828 &$x=2\pi/5,y=26\pi/31,z=58\pi/31$\\
  \hline
  156 & 0.3824 &$x=5\pi/39,y=8\pi/39,z=15\pi/78$\\
  \hline
  158 & 0.3823 &$x=58\pi/79,y=81\pi/79,z=64\pi/79$\\
  \hline
  159 & 0.3814 &$x=8\pi/159,y=64\pi/159,z=30\pi/159$\\
  \hline
  160 & 0.3802 &$x=69\pi/80,y=59\pi/80,z=37\pi/20$\\
  \hline
  162 & 0.3770 &$x=53\pi/21,y=10\pi/9,z=19\pi/81$\\
  \hline
  165 & 0.3760 &$x=24\pi/165,y=26\pi/165,z=34\pi/165$\\
  \hline
  166 & 0.3699 &$x=14\pi/83,y=21\pi/83,z=10\pi/83$\\
  \hline
  169 & 0.3696 &$x=56\pi/169,y=76\pi/169,z=284\pi/169$\\
  \hline
  171 & 0.3678 &$x=32\pi/171,y=294\pi/171,z=6\pi/171$\\
  \hline
  178 & 0.3664 &$x=145\pi/89,y=26\pi/89,z=10\pi/89$\\
  \hline
  180 & 0.3636 &$x=\pi/9,y=97\pi/90,z=127\pi/90$\\
  \hline
  193 & 0.3598 &$x=90\pi/193,y=98\pi/193,z=26\pi/193$\\
  \hline
  204 & 0.3566 &$x=13\pi/51,y=4\pi/51,z=5\pi/34$\\
  \hline
  208 & 0.3501 &$x=\pi/13,y=8\pi/13,z=65\pi/104$\\
  \hline
  214 & 0.3476 &$x=98\pi/107,y=67\pi/107,z=59\pi/107$\\
  \hline
  220 & 0.3459 &$x=19\pi/11,y=163\pi/110,z=121\pi/110$\\
  \hline
  222 & 0.3438 &$x=19\pi/111,y=22\pi/111,z=15\pi/111$\\
  \hline
  225 & 0.3420 &$x=2\pi/225,y=52\pi/225,z=414\pi/225$\\
  \hline
  234 & 0.3410 &$x=4\pi/117,y=24\pi/117,z=43\pi/117$\\
  \hline
  240 & 0.3371 &$x=71\pi/120,y=11\pi/10,z=187\pi/120$\\
  \hline
  244 & 0.3335 &$x=39\pi/122,y=14\pi/61,z=20\pi/61$\\
  \hline
  245 & 0.3305 &$x=16\pi/245,y=186\pi/245,z=46\pi/245$\\
  \hline
  248 & 0.3291 &$x=103\pi/124,y=39\pi/31,z=179\pi/124$\\
  \hline
  259 & 0.3288 &$x=30\pi/259,y=44\pi/259,z=42\pi/259$\\
  \hline
  262 & 0.3274 &$x=142\pi/131,y=215\pi/131,z=87\pi/131$\\
  \hline
  264 & 0.3247 &$x=79\pi/66,y=129\pi/66,z=215\pi/132$\\
  \hline
  276 & 0.3237 &$x=23\pi/138,y=15\pi/69,z=6\pi/69$\\
  \hline
  287 & 0.3188 &$x=6\pi/287,y=76\pi/287,z=28\pi/287$\\
  \hline
  292 & 0.3164 &$x=65\pi/146,y=14\pi/73,z=82\pi/73$\\
  \hline
  295 & 0.3147 &$x=\pi/5,y=50\pi/59,z=22\pi/59$\\
  \hline
  300 & 0.3126 &$x=\pi/75,y=17\pi/150,z=9\pi/25$\\
  \hline
\end{tabular}
\end{center}

\end{ta}

%%%%%%%%%%%%%%%%%%%%%%%%%%%%%%%%%%%%%%%%%%%%%%%%%%%%%%%%%%%%%%%%
\subsection{General Form Constellation Numerical Design}
The connection between the complex Stiefel manifold and $U(M)$
(see the beginning of this section)
makes clear that the techniques used above for square unitary constellations can be
applied to design general form unitary constellations too. For
simplicity we describe the idea with assumption $T=2M$ and
consider the following structure:
$$
\{A^kB|A \in U(T), B=\left(\begin{array}{c}
                  I_M\\
                  0\\
       \end{array}\right), k=0,1, \cdots, L-1\}.
$$
One can check at most $2L-1$ distance calculations are needed to
derive the diversity product (sum or function) with this algebraic
structure.

\begin{ta} The following tables show the constellations (M=2) found
  using SA. More results can be found in~\cite{ha03u2}.

\begin{center}
\begin{tabular}{|c|c|c|c|c|c|c|c|}
  \hline
  size & 3 & 4 & 5 & 6 & 7 & 8 & 9 \\
%  \hline
%  rate & 0.3962& 0.5000 & 0.5905 & 0.6462& 0.7018 & 0.7500 & 0.7925 \\
  \hline
  diversity sum & 0.8654 & 0.7901 & 0.7889 & 0.7652& 0.7514 & 0.7422 & 0.7369 \\
  \hline
  diversity product & 0.8582 & 0.7424 & 0.7330 & 0.6450 & 0.6361 & 0.6216 & 0.5822 \\
  \hline
\end{tabular}

\end{center}
\end{ta}

%%%%%%%%%%%%%%%%%%%%%%%%%%%%%%%%%%%%%%%%%%%%%%%%%%%%%%%%%%%%%%%%%%%%%
\Section{Fast Decoding of the Structured Constellation}
\label{sphere-decoding}

The complexity of ML decoding for unitary space time constellations
increases exponentially with the number of antennas or the
transmission rate. This will preclude its practical use for high
transmission rates or for large number of antennas. Basically our
structured constellations can convert the ML decoding to lattice
decoding naturally, consequently they admit fast decoding
algorithms.

The principle of sphere decoding~\cite{fi85} is as follows:
instead of doing an exhaustive search over all the lattice points,
one can limit its search area to a sphere with given radius
$\sqrt{C}$ centered at received point. One can check the
complexity of this approach in~\cite{fi85} and in~\cite{ha03u1}.

We will use the $A^kB^l$ structure to describe how one can apply
sphere decoding algorithm for the demodulation based on our
constellations. Suppose $A$ has Schur decomposition
$A=U\diag(e^{i\alpha_1}, e^{i\alpha_2},\cdots,e^{i\alpha_M})U^*$,
similarly assume $B=B\diag(e^{i\beta_1},
e^{i\beta_2},\cdots,e^{i\beta_M})B^*$. Consider unitary differential
modulation~\cite{ho00} and denote with $X_{\tau}$ the received
signal at time block $\tau$. The ML demodulation algorithm
involves the following minimization problem:
$$
(\hat{k},\hat{l})=\arg \min_{k,l} {\|X_{\tau}-A^kB^l
X_{\tau-1}\|}_F.
$$
Algebraically one can check that
$$
{\|X_{\tau}-A^kB^l X_{\tau-1}\|}_F={\|A^{-k}X_{\tau}-B^l
X_{\tau-1}\|}_F
$$
$$={\|U\diag(e^{-i k\alpha_1},e^{-i k\alpha_2},
\cdots, e^{-i k\alpha_M})U^*X_{\tau}-V\diag(e^{-i l\beta_1},e^{-i
l\beta_2}, \cdots, e^{-i l\beta_M})V^* X_{\tau-1}\|}_F
$$
So every entry of $X_{\tau}-A^kB^l X_{\tau-1}$ is a
linear combination of trigonometric functions $\cos$ or $\sin$ in
the variables $k, l$, which can be viewed as lattice points.
As demostrated in~\cite{ji03} and ~\cite{ha03u1},
the whole demodulation task has been converted to least-squares problem.
Consequently our structured constellation will admit sphere
decoding algorithm. In~\cite{ji03} a detailed study of the sphere
decoding algorithm applied to constellations from $Sp(2)$ was undertaken.

The complexity (either upper bound or average complexity) of
sphere decoding will depend on the dimension of the lattice. This
will make the weak group structure $A^kB^k$ more remarkable,
because in this case the algorithm requires considering finding the closest
point in a one dimensional lattice, which is very simple.

In~\cite{cl01} a very interesting fast demodulation approach is
proposed for diagonal space time constellations. The authors use
numerical approximation and LLL basis reduction technique to
reduce the decoding complexity. Note that a constellation with the weak group structure
$A^k$ essentially is a diagonal constellation (straightforward Schur
decomposition will show this), therefore the same technique can be
applied to this structure. Most importantly some other algebraic
structure can employ the techniques too. For instance, consider
the $A^kB^lC^m$ structure. If we let $l$ go over a large interval and
let $k, m$ stay within a small interval, the structure will become
``almost'' diagonal. For efficient decoding, one only has to do
exhaustive search for $k, m$ and apply the techniques for diagonal
constellations to decode $l$. Although the decoding complexity
will increase a little, our experiments show the performance will
output the diagonal one remarkably. Exactly the same ``almost''
diagonal idea can be applied to other proposed structures.

%%%%%%%%%%%%%%%%%%%%%%%%%%%%%%%%%%%%%%%%%%%%%%%%%%%%%%%%%%%%
\Section{Conclusions and Future Work}

In this paper, we study the limiting behavior of the {\em
diversity function} by either letting the SNR go to infinity or to
zero. Respectively the {\em diversity product} and the
{\em diversity sum} for unitary constellations are studied from
the analysis of the limiting behavior. We propose algebraic
structures, which are suitable for constructing unitary space
time constellation and feature fast decoding algorithms.
 Based on the presented structure we construct unitary
constellations using geometrical symmetry and numerical methods.
For $2$ dimension most of our codes are better or equal to
the currently existing ones. For higher dimensions many codes
with excellent diversity are found, which were never found before.
Combined with the proposed algebraic structure the numerical
methods can also be employed to optimize the diversity function
at a certain SNR. Future work may involve analyzing the geometric
aspects (such as geodesics, gradients and Hessians of the
functions, etc) on $U(M)$ or the complex Stiefel manifold. Using the
optimization techniques on Riemmannian manifold to optimize the
distance spectrum of a unitary constellation to further search
good-performing constellations is under close investigation too.

%%%%%%%%%%%%%%%%%%%%%%%%%%%%%%%%%%%%%%%%%%%%%%%%%%%%%%%%%%%%%%%%%%%%%%%%

%\nocite{te99,ha97a,sh01,sh02,ha02a,sh00u2,ta00a,ho00a,ho00,li02}
%\nocite{ha02p,ha02p2,aa89,la87b,ph00,no99,ot89,fi85,ha03u1,hu74,ha03u2,bo86b,go98,ea83,ha77}

%\bibliography{huge} \bibliographystyle{plain}

\end{document}